\documentclass[10pt]{amsart}
\usepackage{amsmath,amsthm,amsxtra,amssymb,latexsym,mathrsfs}
\usepackage{verbatim}
\usepackage{amscd}
\usepackage[all]{xy}
\usepackage{eucal}
\usepackage{hyperref}
\usepackage{color}
\usepackage{float}
\usepackage{tikz}
\usepackage[margin=1.5in]{geometry}

\theoremstyle{plain}
\newtheorem{lemma}{Lemma}[section]
\newtheorem{prop}[lemma]{Proposition}
\newtheorem{corollary}[lemma]{Corollary}
\newtheorem{theorem}[lemma]{Theorem}

\newtheorem{THM}{Theorem}
\newtheorem{CORO}[THM]{Corollary}

\theoremstyle{definition}
\newtheorem{definition}[lemma]{Definition}

\newtheorem{remark}[lemma]{Remark}

\newtheorem{example}[lemma]{Example}

\newtheorem*{REM}{Remark}
\newtheorem{DEF}[THM]{Definition}
\newtheorem*{Acknowledgement}{Acknowledgement}

\def \X{\mathfrak{X}}
\def \A{\mathbb{A}}
\def \R{\mathbb{R}}

\def \Z{\mathbb{Z}}
\def \Q{\mathbb{Q}}
\def \N{\mathbb{N}}
\def \G{\mathbb{G}}
\def \P{\mathbb{P}}
\def \p{\mathfrak{p}}
\def \O{\mathcal{O}}
\def \F{\mathcal{F}}
\def \B{\mathcal{B}}
\def \v{\mathfrak{v}}
\def \us{\textup{us}}
\def \ss{\textup{ss}}

\def \In{\textup{in}}
\def \q{\mathfrak{q}}
\def \k{{\bf k}}
\def \r{{\bf r}}

\def \gr{\textup{gr}}
\def \Proj{\textup{Proj}}
\def \Spec{\textup{Spec}}
\def \rank{\textup{rank}}
\def \ord{\textup{ord}}

\begin{document}
\title{On degenerations of projective varieties to complexity-one $T$-varieties}
\author{Kiumars Kaveh}
\address{Department of Mathematics, University of Pittsburgh,
Pittsburgh, PA, USA}
\email{kaveh@pitt.edu}

\author{Christopher Manon}
\address{Department of Mathematics, University of Kentucky, Lexington, KY, USA}
\email{Christopher.Manon@uky.edu}

\author{Takuya Murata}
\address{Department of Mathematics, University of Pittsburgh,
Pittsburgh, PA, USA}
\email{takusi@gmail.com}

\date{\today}

\thanks{The first author is partially supported by a National Science Foundation Grant 
(Grant ID: DMS-1601303), Simons Foundation Collaboration Grant for Mathematicians, and Simons Fellowship.}

\thanks{The second author is partially supported by a National Science Foundation Grant (Grant ID: DMS-1500966).}



\begin{abstract}
Let $R$ be a positively graded finitely generated $\k$-domain with Krull dimension $d+1$. We show that there is a homogeneous valuation $\v: R \setminus \{0\} \to \Z^d$ of rank $d$ such that the associated graded $\gr_\v(R)$ is finitely generated. This then implies that any polarized $d$-dimensional projective variety $X$ has a flat deformation over $\A^1$, with reduced and irreducible fibers, to a polarized projective complexity-one $T$-variety (i.e. a variety with a faithful action of a $(d-1)$-dimensional torus $T$). As an application we conclude that any $d$-dimensional complex smooth projective variety $X$ equipped with an integral K\"ahler form has a proper $(d-1)$-dimensional Hamiltonian torus action on an open dense subset that extends continuously to all of $X$.
\end{abstract}

\maketitle

\setcounter{tocdepth}{1}

\tableofcontents

\section*{Introduction}
Throughout $\k$ denotes an algebraically closed field of characteristic $0$. In this paper we show that if $R$ is a finitely generated positively graded domain of Krull dimension $d+1$ then there is a homogeneous valuation $\v: R \setminus \{0\} \to \Z^d$ of rank $d$ (that is, one less than the Krull dimension) such that the corresponding associated graded $$\gr_\v(R) = \bigoplus_{a \in \Z^d} R_{\v \geq a} / R_{\v > a}$$ is finitely generated (Theorem \ref{th-intro-main} below). In other words, $(R, \v)$ has a finite \emph{Khovanskii basis} (in the sense of \cite[Section 2]{Kaveh-Manon-NOK}). 

Taking $\Proj$, the above then implies the following: Let $L$ be \emph{any} ample line bundle on a projective variety $X$. Then, after going to a sufficiently high tensor power if needed, $L$ can be deformed, in a flat family with reduced and irreducible fibers, to an equivariant ample line bundle $L'$ on a projective complexity-one $T$-variety $X'$ (see Definition \ref{def-intro-flat-degen} and Corollary \ref{cor-intro-toric-degen} below). We recall that a $d$-dimensional variety is a {\it complexity-one} $T$-variety if it possesses a faithful action of the algebraic torus $T=\G_m^{d-1}$ (notice that we do not require the variety to be normal). 

From the results in \cite{HK} one then concludes that any complex $d$-dimensional smooth projective variety with an integral K\"ahler form admits a faithful $(d-1)$-dimensional Hamiltonian torus action on a dense open subset that extends continuously to the whole variety (Theorem \ref{th-intro-Hamiltonian} below). 
 
It is well-known that if $\X \subset \P^N \times \A^1 \to \A^1$ is a flat family of schemes embedded in a projective space $\P^N$ then the Hilbert polynomials of fibers are all the same and thus, in particular, all fibers have the same degree. This illustrates why flat deformations are such useful tools in algebraic geometry, specially when dealing with intersection theoretic data. In our setting, the important and key requirement that fibers are reduced guarantees that we do not need to worry about multiplicities and embedded components at the special fiber. Whenever there is a degeneration of a variety to a variety with a torus action one can gain insight into the geometry of the variety using the torus action. The largest torus action is when the special fiber is a toric variety. Such a degeneration is usually called a \emph{toric degeneration}.
%The larger the torus action, the more structure that we can exploit (one example is \cite{HK} mentioned above and Theorem \ref{th-intro-Hamiltonian} later in the introduction). The largest torus action is when the special fiber is a toric variety. Such a degeneration is usually called a \emph{toric degeneration}.
%Naturally it is more desirable to have a degeneration of a given projective variety to a projective toric variety rather than a complexity-one $T$-variety (namely a {\it toric degeneration}). 
Toric degenerations of different special classes of varieties have been constructed and studied by many authors during the past few decades. See for example \cite{FFL} for a nice survey of toric degenerations of flag varieties and connections with representation theory, and \cite{Alexeev-Brion, Kaveh-Crystal} for toric degenerations of spherical varieties. 

It is an important and fundamental question \emph{whether an ample line bundle on a projective variety has a flat deformation, with reduced and irreducible fibers, to an equivariant ample line bundle on a (not necessarily normal) projective  toric variety} (i.e. a toric degeneration). 
\begin{itemize}
\item The main result of the present paper (Theorem \ref{th-intro-main} and Corollary \ref{cor-intro-toric-degen}) shows that the answer to this question is \emph{almost} affirmative in the sense that we can always degenerate to a complexity-one $T$-variety. %Moreover, we give examples of embeddings of projective curves that do not have \emph{any} toric degenerations.
\item On the other hand, we give examples of line bundles on projective curves that do not have any toric degenerations. More precisely we show the following: Let $D$ be a very ample divisor on a smooth complex projective curve $X$ with genus $> 1$. Then the corresponding line bundle $L=\mathcal{O}(D)$ has a toric degeneration 
%the algebra of sections of $\mathcal{O}(D)$ has a homogeneous valuation with finitely generated value semigroup 
if and only if some multiple of $D$ is equivalent to a multiple of a point (see Section \ref{sec-curve-case} and Propositions \ref{prop-curve-finitely-gen} and \ref{prop-curve-not-finitely-gen}). %This example seems to have been known to Dave Anderson as well.
\end{itemize}

In addition, we show that any toric degeneration (in the sense of Definition \ref{def-intro-flat-degen} below) comes from a valuation with a finite Khovanskii basis (Theorem \ref{th-degen-valuation}). %We should point out that it is possible to find rational curves with a projective coordinate ring which cannot have a toric degeneration.
In regard to examples, in \cite[Corollary 3.14]{Ilten-Wrobel} it is also shown that the projective coordinate ring of a very general integral rational plane curve of degree greater than $3$ has no homogeneous valuation with a finite Khovanskii basis and hence has no toric degeneration by our Theorem \ref{th-degen-valuation}. 

Toric varieties can be classified using combinatorial and convex geometric data such as lattice points and convex polytopes. In the same spirit, there is a classification of normal complexity-one $T$-varieties using the so-called polyhedral divisors (see \cite{Altman}). Normal complexity-one $T$-varieties over a rational curve are particularly nice and admit many degenerations to toric varieties (see \cite{Ilten-Manon}).

Let us start by giving a precise definition of the notion of degeneration that we use throughout the paper. %Let $\k$ be the base field which we assume to be algebraically closed. Let $X$ be a variety of dimension $d$ over $\k$. 
Let $R$, $R'$ be finitely generated positively graded $\k$-domains with Krull dimension $d+1$. Let $X = \Proj(R)$, $X' = \Proj(R')$ be the corresponding $d$-dimensional projective varieties equipped with sheaves of modules $\O_X(n)$ and $\O_{X'}(n)$, $n \in \N$, respectively. 
\begin{DEF}   \label{def-intro-flat-degen}
By a {\it flat degeneration} (or \textit{degeneration} for short) of $X$ to $X'$ we mean a $\G_m$-equivariant flat family $\pi: \X \to \A^1$ of varieties with reduced fibers and equipped with sheaves of $\O_\X$-modules $\O_\X(n)$ such that:
\begin{itemize}
%\item[(a)] The $\G_m$-action on $\X$ that lifts the standard action of $\G_m$ on $\A^1$. This in particular implies that the family is trivial outside $X_0 = \pi^{-1}(0)$, that is, $\pi: \X \setminus \pi^{-1}(0) \to \A^1 \setminus \{0\}$ is trivial with fiber $X$. Hence for $t \neq 0$, the fibers $X_t = \pi^{-1}(t)$ are all isomorphic to $X$.
\item[(a)] The fibers $X_t := \pi^{-1}(t)$, $t \neq 0$, are all isomorphic to $X$ and the special fiber $X_0 := \pi^{-1}(0)$ is isomorphic to $X'$. Note that the $\G_m$-equivariance of the family then implies that $\X \setminus X_0$ is a trivial family with fiber $X$.
\item[(b)] Under the isomorphism $X_1 \cong X$, the sheaf $\O_\X(n)_{|X_1}$ coincides with $\O_X(n)$, for all $n \in \N$.
\item[(c)] Under the isomorphism $X_0 \cong X'$, the sheaf $\O_\X(n)_{|X_0}$ coincides with $\O_{X'}(n)$, for all $n \in \N$.
 %\item[(a)] The special fiber $X_0 = \pi^{-1}(0)$ is isomorphic to $X'$.
%\item[(b)] $\O_\X(1)$ restricted to $X_1$ (respectively $X_0$) is isomorphic to $\mathcal{O}_X(1)$ on $X$ (respectively $\O_{X'}(1)$ on $X'$).
\end{itemize}
%The $\G_m$-equivariance of the family in particular implies that the family is trivial outside $X_0$, that is, $\pi: \X \setminus \pi^{-1}(0) \to \A^1 \setminus \{0\}$ is trivial with fiber $X$. Hence for $t \neq 0$, the fibers $X_t = \pi^{-1}(t)$ are all isomorphic to $X$.
We write $X \rightsquigarrow X'$ when there is a degeneration of $X$ to $X'$. We call a degeneration as above a {\it toric degeneration of $X$} if in addition $X'$ is a (not necessarily normal) toric variety and the $\O_{X'}(n)$ are torus equivariant.
\end{DEF}

%We should point out that in this paper, we do not assume normality in the definition of a toric variety. Thus for us a toric variety is a variety that has a torus action with an open orbit.

%{\bf CHRIS: Should we include the condition that the $\mathbb{G}_m$ action on $\mathbb{A}^1$ lifts to the total family of the degeneration in the geometric description above?  Note that this is implicit in the algebraic description.}

We remark that our definition of a flat degeneration above is very close to the notion of \emph{test configuration} from K\"ahler geometry. 

Let $R$ be a finitely generated positively graded algebra and domain with $\dim(R)=d+1$ and $R_0 = \k$.
%Let $X = \Proj(R)$ and $d=\dim(X)$. 
The main result of the paper which enables us to construct degenerations to complexity-one $T$-varieties is the following (Theorem \ref{th-main}):
\begin{THM}   \label{th-intro-main}
There exists a homogeneous rank $d$ valuation $\mathfrak{v}: R \setminus\{0\} \to \Z_{\geq 0}  \times \Z^{d-1}$ such that the associated graded $\gr_\v(R)$ is finitely generated (in other words, $(R, \mathfrak{v})$ has a finite Khovanskii basis). Here the first component of $\mathfrak{v}$ is the degree with respect to the grading of $R$, also $\Z_{\geq 0} \times \Z^{d-1}$ is ordered lexicographically and the ordering on $\Z_{\geq 0}$ is the reverse of its natural ordering.
\end{THM}

Taking $\Proj$ and in light of \cite[Proposition 5.1]{Anderson} we obtain the following:
\begin{CORO}   \label{cor-intro-toric-degen}
Let $X \subset \P^N$ be a projective variety with homogeneous coordinate ring $R = \k[X]$. 
%Let $X$ be a projective variety with an embedding $X \hookrightarrow \P^N$ and $R = \k[X]$ the corresponding homogeneous coordinate ring. 
Then $X$ can be degenerated to a complexity-one $T$-variety $X'$ via a $\G_m$-equivariant family $\pi: \X \to \A^1$ (as in Definition \ref{def-intro-flat-degen}). Moreover, the family $\X$ can be embedded in $W\P \times \A^1$ where $W\P$ is a weighted projective space and the following holds:
\begin{itemize}
\item[(1)] The weighted projective space $W\P  = \P(1, \ldots, 1, \lambda_1, \ldots, \lambda_s)$ contains $\P^N = \P(1, \ldots, 1)$.
\item[(2)] $\pi: \X \to \A^1$ is given by the projection on the second factor $W\P \times \A^1 \to \A^1$.
\item[(3)] The $\G_m$-action on $\X$ is the restriction of a linear $\G_m$-action on $W\P \times \A^1$.
\item[(4)] The sheaves $\O_\X(n)$ are all induced from the Serre sheaves $\O_{W\P}(n)$ on the weighted projective space via the embedding $\X \hookrightarrow W\P \times \A^1$. 
\item[(5)] The $T$-action on $X_0 \cong X'$ comes from a linear $T$-action on $W\P = W\P \times \{0\}$. 
\end{itemize}

%{\bf Equivalently, there is a projective space $\P$ with a very ample line bundle $\mathcal{L}$ such that $L_0 = \mathcal{L}_{|X_0}$ is $T$-linearlized and $
%\mathcal{L}_{|X_1}$ is a tensor power of the line bundle $L$ on $X$.}  

%Let $X$ be a $d$-dimensional projective variety $X$. Then 
%any embedding $X \hookrightarrow \P^N$ can be degenerated to an equivariant embedding of a complexity-one $T$-variety $X'$ into a weighted projective space $W\P$ where $T$ acts on $W\P$ via a linear action. More specifically, the weighted projective space $W\P  = \P(1, \ldots, 1, \lambda_1, \ldots, \lambda_s)$ contains $\P^N = \P(1, \ldots, 1)$ and we have the following: 
%\begin{itemize}
%\item[(a)] There is a $\G_m$-equivariant flat family $\pi: \X \to \A^1$ with reduced fibers where $\X \subset W\P \times \A^1$ and $\pi$ is induced by the projection on the second factor $W\P \times \A^1 \to \A^1$. 
%\item[(b)] The embedding $X_1 \hookrightarrow W\P \times \{1\}$ coincides with $X \hookrightarrow \P^N \subset W\P$.
%\item[(c)] The embedding $X' \hookrightarrow W\P$ coincides with $X_0 \hookrightarrow W\P \times \{0\}$.
%\end{itemize}
\end{CORO}

%\begin{REM}
%The first version of the paper claimed that any polarized projective variety can be degenerated in a sequence of degenerations to an equivariantly polarized projective toric variety. As mentioned above (see Example \ref{prop-curve-not-finitely-gen}) this claim is false. 
%%Using our methods we can prove only degeneration to a complexity-one $T$-variety. 
%Our method of constructing degenerations in this paper essentially relies on the Bertini irreducibility theorem (\cite[Theorem 0.1]{Bertini}) which fails in dimension $1$ and hence in general we can achieve only degeneration to a complexity-one $T$-variety. 
%\end{REM}

\begin{REM}
It seems that the proofs of Theorem \ref{th-intro-main} and Corollary \ref{cor-intro-toric-degen} work over a characteristic $p$ base field as well but the authors have not checked all the details. 
\end{REM}

\begin{REM}
If $R$ is Cohen-Macaulay, there is a relatively simple construction of a complexity-one degeneration of $X = \Proj(R)$. If moreover, $X$ is assumed to be smooth and the genus of the curve $C = V(f_1, \ldots, f_{d-1})$, where the $f_i \in R_1$ are in general position, is zero or one then we can complete our complexity-one degeneration to a toric degeneration (see Section \ref{sec-CM-case}).
\end{REM}

Below we briefly mention some applications of our main result.

\subsubsection*{Hamiltonian torus actions} Combining the main result in \cite{HK} and Corollary \ref{cor-intro-toric-degen} we conclude the following general result about constructing Hamiltonian torus actions on arbitrary smooth projective varieties.

\begin{THM}  \label{th-intro-Hamiltonian}
Every smooth projective variety of dimension $d$ equipped with an integral K\"ahler form has a proper faithful $(d-1)$-dimensional Hamiltonian torus action on a dense open subset that extends continuously to the whole variety. Moreover, the image of the moment map is the standard simplex. 
\end{THM}

In a series of papers, Karshon and Tolman classify proper complexity-one Hamiltonian torus actions (see \cite{KT}). In their terminology, the Hamiltonian action in Theorem \ref{th-intro-Hamiltonian} is {\it tall}. 

%Degenerations to $T$-varieties and in particular toric degenerations can also be used to approach questions in geometric quantization theory. 
Corollary \ref{cor-intro-toric-degen} is expected to have applications in proving general results in geometric quantization theory. To this end one needs to extend some known results about toric varieties to complexity-one $T$-varieties. More specifically, building on the work \cite{HK}, in \cite{HHK} it is shown that given a smooth projective variety together with a quantization data as well as a toric degeneration, the ``real polarization'' converges to the ``K\"ahler polarization'', in a sense that is made precise in there.

\subsubsection*{Tropical geometry and existence of prime cones}
Let $A$ be a finitely generated positively graded $\k$-domain presented as $A \cong \k[x_1, \ldots, x_n] / I$ where $I$ is a homogeneous ideal. It is well-known that tropical variety $\text{trop}(I)$ of the ideal $I$ is the support of a fan in $\R^n$ of pure dimension $d = \dim(A)$. A cone in such a fan is called a \emph{prime cone} if the corresponding initial ideal of $I$ is a prime ideal. 
It is a problem of great importance to understand prime cones in a presentation of $A$. In particular one would like to know whether $A$ has a presentation such that the corresponding tropical variety has a prime cone, and more strongly a prime cone of maximal dimension. This is already a very difficult question even in specific examples such as the Pl\"ucker algebra of Grassmannian $\text{Gr}(k, N)$ for $k \geq 3$ (see \cite{Mohammadi-Shaw} as well as \cite{BLMM}). It is known that maximal cones in the tropical variety of $\text{Gr}(2, N)$ are all prime. Prime cones give rise to degenerations of $\Proj(A)$ by degenerating $I$ to its initial ideal. In particular, a maximal prime cone gives a toric degeneration.

The main result of the present paper combined with \cite[Theorem 2]{Kaveh-Manon-NOK} gives the following strong existence result:
\begin{THM}
Let $\{f_1, \ldots, f_r\}$ be a set of algebra generators for a $\k$-domain $A$. Then this set can be enlarged to $\{f_1, \ldots, f_s\}$ such that the ideal $I$ of relations among the $f_i$ has a prime cone of codimension at most $1$ (i.e. a prime cone of dimension $d-1$ or $d$).  
\end{THM}

\subsubsection*{Computational algebra} 
The notion of a Khovanskii basis far generalizes that of a SAGBI basis (Subalgebra Analogue of Gr\"obner Basis for Ideals) for subalgebras of a polynomial algebra to arbitrary domains equipped with a valuation \cite[Section 2]{Kaveh-Manon-NOK}. This theory allows to extend the scope of Gr\"obner basis methods for polynomial rings to a far larger class of algebras. While there are examples of subalgebras that do not have a finite SAGBI basis, our Theorem \ref{th-intro-main} shows that \emph{any} positively graded finitely generated domain has a valuation with rank at most one less than the Krull dimension which possesses a finite Khovanskii basis.

\begin{REM} 
The general theory of Newton-Okounkov bodies associates a convex body $\Delta \subset \R^d$ (called a {\it Newton-Okounkov body}) to the algebra of sections $R$ of a line bundle $L$ on a $d$-dimensional projective variety $X$ (see \cite{KKh, LM}). More generally, this construction associates convex bodies to a positively graded algebra $R$. The dimension and volume of the convex body $\Delta$ encode asymptotic information about the Hilbert function of $R$ (see \cite{KKh}). The construction of a Newton-Okounkov body $\Delta$ requires the extra choice of a full rank valuation $\mathfrak{v}': \k(X) \setminus \{0\} \to \Z^d$. The valuation $\mathfrak{v}'$ can be naturally extended to a valuation $\mathfrak{v}: R \setminus \{0\} \to \Z_{\geq 0} \times \Z^d$ where the first component is given by the grading of $R$. 
The convex body $\Delta$ is constructed from the value semigroup 
$S=S(R, \mathfrak{v})$ of values of $\mathfrak{v}$ on $R \setminus \{0\}$. When this semigroup is finitely generated, the convex body $\Delta$ is a rational convex polytope. 
%Given an algebra $R$ with a valuation $\v$ the corresponding associated graded $\gr_\v(R) = \bigoplus_{a} R_{\v \geq a} / R_{\v > a}$ is by construction graded by the value semigroup $S(R, \v)$ and hence in general is more structured than the original algebra $R$. If the valuation $\v$ has full rank (that is, its rank is equal to the Krull dimension of $R$) and the base field $\k$ is algebraically closed, the associated graded $\gr_\v(R)$ is isomorphic (as a graded algebra) to the semigroup algebra $\k[S(R, \v)]$. 
\end{REM}

%\begin{REM}
%{Deformations of algebraic varieties and in particular deformations of curves are well-studied, in particular in the context of test configurations. We can mention \cite{Pinkham} in this regard. Nevertheless we are not aware of any general construction of a complexity-one degeneration in the literature as in Corollary \ref{cor-intro-toric-degen}.} 
%In particular, some of the material about curves in Section \ref{sec-curve-case} might be known to the experts. 
%\end{REM}

\begin{REM}
We would like to mention \cite[Proposition 14]{AKL} which gives a sufficient condition for the existence of a full rank valuation with a finitely generated value semigroup for the algebra of sections of a very ample line bundle on a smooth projective variety (which then implies the existence of a toric degeneration of this variety). This is related to our Proposition \ref{prop-curve-finitely-gen}. As they point out this  condition is far from necessary and often does not hold. 
%Similar to our approach, their method also relies on a Bertini type argument although for the purposes of 
Unlike the present paper, \cite{AKL} does not consider non-full rank valuations or non-full dimensional torus actions on the special fiber. For a non-full rank valuation the finite generation of the associated graded algebra is a much stronger condition than the finite generation of the value semigroup and is a source of technicalities appearing in the proof of our main theorem. Also we point out that in our main result (Corollary \ref{cor-intro-toric-degen}) we do not assume that the variety is smooth.   
\end{REM}

\begin{REM}
The recent work \cite{Ilten-Wrobel} among other things gives tools for checking if a complexity-one variety has a toric degeneration. In \cite[Theorem 1.1]{Ilten-Wrobel}, they show that verifying this property for a given valuation on a complexity-one variety essentially reduces to a check in the curve case. Moreover, if the complexity-one degeneration is rational, \cite[Theorem 4.1]{Ilten-Wrobel} says one must only check this property holds for a finite number of rational curves.  For this latter case see also \cite[Theorem 3.5]{Ilten-Wrobel}. 
\end{REM}

%Whenever the associated graded $\gr_\v(R)$ is a finitely generated algebra one says that $(R, \v)$ has a finite {\it Khovanskii basis} (see Section \ref{subsec-Khovanskii-bases} and \cite{Kaveh-Manon-NOK}). 
%This is far extension of notion of SAGBI basis for subalgebras of a  polynomial ring (see \cite{Robbiano-Sweedler}).  
%If the valuation is full rank and the base field is algebraically closed this is equivalent to $S(R, \mathfrak{v})$ be a  finitely generated semigroup. In principle, when we have a finite  Khovanskii basis, we can do computations in the algebra more effectively and algorithmically (see \cite[Section 2]{Kaveh-Manon-NOK}). We would like to point out that, given a valuation on an algebra, firstly the value semigroup is often not finitely generated and secondly it is usually not easy to verify the finite generation of the value semigroup.

Finally we discuss the organization of the paper and outline of the proof of the main results. Section \ref{sec-prelim} covers some needed background material on filtrations, valuations, Rees algebras and symbolic powers of ideals. Section \ref{subsec-Khovanskii-bases} talks about Khovanskii bases and degenerations coming from valuations. A rather new result in this section is Theorem \ref{th-degen-valuation} that shows toric degenerations necessarily come from full rank valuations. 

In Section \ref{sec-degen-stages} we prove the main step in the proof of Theorem \ref{th-intro-main}, namely we show that if $X = \Proj(R)$ is a projective variety of dimension $d$ then there is a sequence of degenerations $X_0 = X \rightsquigarrow X_1 \rightsquigarrow \cdots \rightsquigarrow X_{d-1}$ that degenerates $X$ to a complexity-one $T$-variety $X' = X_{d-1}$. We call this {\it degeneration in stages}. The main idea behind degeneration in stages is the Bertini irreducibility theorem. The key technical lemma in the proof of degeneration in stages is a statement about finiteness of symbolic Rees algebras of certain height one prime ideals. This is proved in the appendix (Section \ref{sec-appendix}). 

In Section \ref{sec-curve-case} we consider the case where $X$ is a projective curve. We give examples of line bundles on curves that do not have any toric degenerations. Section \ref{sec-CM-case} discusses the particularly nice case where $R$ is Cohen-Macaulay (or in other words, $X$ is arithmetically Cohen-Macaulay). In this case, by a well-known theorem of Rees about associated graded algebra of an ideal generated by a regular sequence (Theorem \ref{th-Rees}), we can construct a degeneration of $X$ to a certain  complexity-one $T$-variety which is a compactification of a trivial vector bundle over a curve. 

In Section \ref{sec-finite-Khovanskii-basis} we show that a sequence of degenerations (as constructed in Section \ref{sec-degen-stages}) gives rise to a rank $d$ valuation 
on $R$ with finitely generated associated graded algebra (in other words with a finite Khovanskii basis). This in turn gives a degeneration of $X$ to a complexity-one $T$-variety (in one step) proving Theorem \ref{th-intro-main}. 
 
\begin{Acknowledgement}
Some of the main results of this paper are extracted from the third author's  Ph.D. thesis in progress. We would like to thank Craig Huneke, Javid Validashti, Dave Anderson, Mircea Musta\c{t}{a}, Dale Cutkosky, Michel Brion, Igor Dolgachev, Yael Karshon, Maksym Fedorchuk, Bernd Sturmfels and Mateusz Michalek for helpful discussions and correspondence. In particular, we thank Craig Huneke for pointing out a shorter proof for Lemma \ref{lem-main}.
%We are also grateful to Michel Brion for pointing to us that a Hilbert scheme may not contain a toric variety. 
Finally thanks to Dave Anderson, Lara Bossinger and Diane Maclagan for reading a very first draft of this paper. 
\end{Acknowledgement}

\section{Preliminaries on valuations, Rees algebras and symbolic powers} \label{sec-prelim}
Let $R$ be a finitely generated $\k$-domain. In this section we discuss the Rees algebra and associated graded corresponding to a filtration, and in particular to a valuation, on $R$.
%\footnote{We point out that we only use valuations with values in $\Z$.} 
%We prove a general statement for when these algebras are finitely generated (Theorem \ref{th-finite-generation-symbolic-gr}). When the Rees algebra is finitely generated we have a deformation of our original algebra $R$ to the associated graded. %We refer to $\Proj$ of the symbolic associated graded as the {\it symbolic normal cone}. 
%This is an important step in the main construction of the paper (Theorem 
%\ref{th-degeneration-in-stages}). 
\subsection{Generalities on multiplicative filtrations} \label{subsec-generalities-filt}
%Let $R$ be a finitely generated algebra and domain. 
A {\it multiplicative filtration} $\F = (F_i)_{i  \in \Z}$ in $R$ is a descending sequence of vector subspaces $\cdots \supset F_{i-1} \supset F_i \supset F_{i+1} \supset \cdots $ in $R$ such that for all $i, j \geq 0$ we have 
%\begin{equation}   \label{equ-multiplicative}
$F_i F_j \subset F_{i+j}$.
%\end{equation}
We also assume that $F_i = R$ for all $i \leq 0$. Given a multiplicative filtration $\F$ on $R$, one defines its corresponding {\it associated graded} $\gr_\F(R)$ by:
$$ \gr_\F(R) = \bigoplus_{i} F_i / F_{i+1}.$$
Also one defines its {\it Rees algebra} to be the algebra: $$\mathcal{A}_\F(R) = \bigoplus_{i} F_i.$$
The multiplicativity assumption on $\F$ guarantees that these are actually algebras. We note that the Rees algebra $\mathcal{A}_\F(R)$ is moreover a $\k[t]$-module, where $\k[t]$ is the polynomial algebra in one indeterminate $t$. The $\k[t]$-module structure is defined as follows: for $f \in F_i$ we let $t \cdot f$ to be $f$ regarded as an element of $F_{i-1}$. %We also define the action of $t$ on $F_0$ to be $0$. 
One thinks of the Rees algebra $\mathcal{A}_\F(R)$ as a family of algebras that deform the algebra $R$ to the associated graded $\gr_\F(R)$ (see \cite[Proposition 2.2]{Teissier}).

Let $0 \neq f \in R$ and let $i \in Z_{\geq 0}$ be such that $f \in F_i$ but $f \notin F_{i+1}$. We denote the image of $f$ in $F_i / F_{i+1} \subset \gr_\F(R)$ by $\bar{f}$. One verifies that $f \mapsto \bar{f}$ is a multiplicative homomorphism from $R \setminus \{0\}$ to the multiplicative set of homogeneous elements of $\gr_\F(R)$.

One can give examples of a filtration $\F$ on a finitely generated algebra $R$ such that the Rees algebra $\mathcal{A}_\F(R)$ and the associated graded $\gr_\F(R)$ are not finitely generated algebras.

A main source of examples of multiplicative filtrations on an algebra are powers of ideals. Let $I \subset R$ be an ideal. It is clear that the sequence of subspaces $F_i = I^i$  (where $I^i = R$ for $i \leq 0$) is a multiplicative filtration $\F$ on $R$. We denote the corresponding associated graded and the Rees algebra by $\gr_I(R)$ and $\mathcal{A}_I(R)$ respectively. It is well-known that these algebras are finitely generated. In algebraic geometry language, $\Spec(\gr_I(R))$ is the normal cone to the subscheme corresponding to $I$ in $\Spec(R)$. The scheme $\Spec(\mathcal{A}_I(R))$ gives the deformation of $\Spec(R)$ into the normal cone $\Spec(\gr_I(R))$ and is related to the notion of blowup along the ideal $I$ (see \cite[Chapter 5]{Fulton}). When $R$ is a positively graded algebra one can take $\Proj$ of these algebras to get the projective version of deformation to normal cone.

While the algebras $\gr_I(R)$ and $\mathcal{A}_I(R)$ are always finitely generated, it is not difficult to come up with an example in which $\gr_I(R)$ is not an integral domain (see Example \ref{elliptic curve}). Since we are interested in degenerating a given variety to another variety with a torus action, we would like to have filtrations for which the associated graded is a domain. For this reason we consider filtrations associated to valuations. 

More generally, one defines a multiplicative filtration indexed by an ordered group. By an ordered group we mean an abelian group $\Gamma$ (written additively) equipped with a total order $>$ which respects the group operation, i.e. for $a, b, c \in \Gamma$ we have $a > b$ implies that $a+c > b+c$. A {\it multiplicative filtration indexed by $\Gamma$} is a collection 
$\F = (F_a)_{a \in \Gamma}$ of vector subspaces of $R$ such that $F_a \subset F_b$ whenever $a > b$ (i.e. the filtration is decreasing). Moreover, for $a, b \in \Gamma$ we have $F_aF_b \subset F_{a+b}$.
We assume that $F_a = R$ for $a < 0$, where $0$ denotes the identity element in $\Gamma$.  

Similar to filtrations indexed by $\Z$, we define the {\it associated graded} of $\F$ by $\gr_\F(R) = \bigoplus_{a \in \Gamma} F_a / F_{> a}$. Here $F_{> a} = \sum_{b>a} F_b$. We also define its {\it Rees algebra} to be $\mathcal{A}_\F(R) = \bigoplus_{a \in \Gamma} F_a$. It follows from the definition that the algebras $\gr_\F(R)$ and $\mathcal{A}_\F(R)$ are graded by the group $\Gamma$.

%Let $\Gamma_{> 0} = \{ a \in \Gamma \mid a > 0\}$. Let us suppose that $\Gamma_{> 0}$ is maximum well-ordered i.e. any increasing chain stops. 
Suppose for any $0 \neq f \in R$ the set $\{b \mid f \in F_b\}$ has a maximum $a$. Then $f \in F_a$ but $f \notin F_{> a}$. Similar to before, we denote the image of $f$ in $F_a / F_{>a} \subset \gr_\F(R)$ by $\bar{f}$. Again $f \mapsto \bar{f}$ is a multiplicative homomorphism from $R \setminus \{0\}$ to the multiplicative set of homogeneous elements of $\gr_\F(R)$.

%the notion of symbolic powers of an ideal (instead of the usual powers of the ideal). 

%Before we discussing filtration associated to a valuation, we give the following criterion for when the associated graded of a filtration is a domain.
%\subsection{Symbolic powers of an ideal}

\subsection{Valuations}   \label{subsec-valuations}
We recall some generalities about valuations and their associated filtrations. Let $R$ be a $\k$-algebra and domain, and let $\Gamma$ be an ordered group. Throughout the paper, the ordered groups we work with are $\Z$ or $\Q$ (with their standard ordering) or $\Z^r$, $\Q^r$ (for some $r > 0$) and equipped with some total order (often a lexicographic order).

\begin{definition}
A function $\v: R \setminus \{0\} \to \Gamma$ is a {\it quasivaluation} on $R$ if the following hold:
\begin{itemize}
\item[(1)] For all $0 \neq f, g, f+g \in R$ we have $\v(f + g) \geq \min\{\v(f), \v(g)\}$.
\item[(2)] For all $0 \neq f, g \in R$ we have $\v(fg) \geq \v(f) + \v(g)$.
\item[(3)] For all $0 \neq f \in R$ and $0 \neq c \in \k$ we have $\v(cf) = \v(f)$.
\end{itemize}
If for all $0 \neq f, g \in R$ we have $\v(fg) = \v(f) + \v(g)$ then $\v$ is called a {\it valuation}.
\end{definition}
It is sometimes useful to define a quasivaluation to be a map $\v: R \to \Gamma \cup \{\infty\}$ satisfying the above axioms, where $\infty$ is greater than all elements in $\Gamma$.

In the case that $\v$ is a valuation, the image of $R\setminus \{0\}$ under $\v$ is a subsemigroup of the value group $\Gamma$; we call this the {\it value semigroup} of $(R, \v)$, and we denote it by $S(R, \v)$. 

In this paper we only work with discrete valuations, that is, the image of the valuation is a discrete subset of $\Q^r$. The {\it rank} (also called the {\it rational rank}) of a valuation is the rank of the subgroup generated by the image of the valuation. One shows that the rank of a valuation $\v: R \setminus \{0\} \to \Q^n$ is less than or equal to $\dim(R)$, the Krull dimension of $R$ (see \cite[Chapter VI]{SZ}). We call $\v$ a {\it full rank valuation} if the rank of $\v$ is equal to $\dim(R)$. %For most of the paper we work with valuations whose value group is $\Z$ (rank $1$ valuations), but in Section \ref{sec-finite-Khovanskii-basis} we will deal with higher rank valuations.  

Let $\F = (F_a)_{a \in \Gamma}$ be a multiplicative filtration. Define the function $\v_\F: R \to \Gamma \cup \{\infty\}$ as follows: for $0 \neq f \in R$,
\begin{equation}  \label{equ-v_F}
\v_\F(f) = \max\{ a \mid f \in F_a \}.
\end{equation}
If the maximum is not attained we let $\v_\F(f) = \infty$, in particular, $\v(0) = \infty$. One verifies that $\v_\F$ is a quasivaluation on the algebra $R$. 
%We recall that a function $v: R \to \Z \cup \{\infty\}$ is a {\it quasivaluation} if it satisfies the following: (1) for all $x, y \in R$ we have $v(x+y) \geq \min\{ v(x), v(y) \}$. (2) For all $x, y \in R$, $v(xy) \geq v(x) + v(y)$. If $v$ instead satisfies the stronger condition that $v(xy) = v(x)+v(y)$, for all $x, y \in R$, then $v$ is a {\it valuation}.
Conversely, if $\v: R \to \Gamma \cup \{\infty\}$ is a quasivaluation, one can define a filtration $\F_\v = (F_{\v \geq a})_{a \in \Gamma}$ as follows. For each $a$, let:
$$F_{\v \geq a} = \{ 0 \neq f \in R \mid \v(f) \geq a \} \cup \{0\}.$$  
One checks that $\F_\v$ is indeed a multiplicative filtration, and moreover, the operations of $\F \mapsto \v_\F$ and $\v \mapsto \F_\v$ are inverse of each other.

\begin{example}   \label{ex-grading-val}
%Suppose $R = \bigoplus_{i \geq 0} R_i$ be a positively graded algebra. For $f = \sum_i f_i$ define $v(f) = \min\{i \mid f_i \neq 0\}$. Then $v: R \setminus \{0\} \to \Z_{\geq 0\}$ is a valuation.
Suppose $R = \bigoplus_{a \in \Gamma} R_a$ is an algebra graded by an ordered abelian group $\Gamma$. Suppose the support semigroup $S(R, \Gamma) = \{ a \in \Gamma \mid R_a \neq \{0\} \}$ is maximum well-ordered. For $f = \sum_a f_a \in R$ define $\v(f) = \min\{ a \mid f_a \neq 0\}$. Then $\v: R \setminus \{0\} \to \Gamma$ is a valuation with value semigroup $S(R, \Gamma)$. A common example of this situation is when $R$ is a positively graded algebra. 
\end{example}

%Definition of valuation homogeneous with respect to a grading
Letting $R = \bigoplus_{a \in \Lambda} R_a$ be an algebra graded by an abelian group $\Lambda$, we can consider any homomorphism $\pi: \Lambda \to \Gamma$ of abelian groups, where $\Gamma$ is an ordered group.  The function $\v_{\pi}: R \setminus \{0\} \to \Gamma$ which sends $f = \sum_a f_a$ to $\min\{ \pi(a) \mid f_a \neq 0\}$ is easily checked to define a valuation on $R$.  We call functions which arise in this way {\it grading functions}.

More generally, let $R = \bigoplus_{a \in \Lambda} R_a$ be an algebra graded by an abelian group $\Lambda$ and let $\v: R \setminus \{0\} \to \Gamma$ be a valuation.  We say that $\v$ is {\it homogeneous with respect to the $\Lambda$-grading} if for any $f = \sum_a f_a \in R$ we have: 
$$\v(f) =  \min\{ \v(f_a) \mid f_a \neq 0\}.$$

Given a quasivaluation $\v$ on $R$ we denote the Rees algebra and the associated graded algebra corresponding to the filtration $\F_\v$ by $\mathcal{A}_\v(R)$ and $\gr_\v(R)$ respectively.  The following propositions are straightforward to verify from the definitions.
\begin{prop}   \label{prop-ass-graded-domain}
Let $\F$ be a multiplicative filtration on a domain $R$ with corresponding quasivaluation $\v_\F$ (given by \eqref{equ-v_F}). Then $\gr_\F(R)$ is a domain if and only if $\v_\F$ is a valuation, that is, for any $0 \neq f, g \in R$ we have $\v_\F(fg) = \v_\F(f) + \v_\F(g)$.
\end{prop}

\begin{prop} \label{prop-ass-graded-graded}
Suppose $R$ is graded by an abelian group $\Lambda$ as above, and $\v$ is homogeneous with respect to this grading.  Then the algebras $\mathcal{A}_\v(R)$ and $\gr_\v(R)$ are graded by $\Lambda$, and the support semigroups (see Example \ref{ex-grading-val}) $S(R, \Lambda)$, $S(\mathcal{A}_\v(R), \Lambda)$ and $S(\gr_\v(R), \Lambda)$ all coincide. Furthermore, if $\v$ is a grading function, the algebras $\gr_\v(R)$ and $R$ are seen to be isomorphic as $\Lambda$-graded algebras. 
\end{prop} 

%{\bf NOTE: For clarity we should emphasize the motivation behind the degeneration-in-stages in this language: at each step we must first preserve the grading that has already been constructed by selecting a homogeneous valuation (support semigroups coincide), and second we must ensure that we add additional grading, ie we must take care that we don't choose a grading function (associated graded is canonically isomorphic to the original algebra).}

\subsection{Symbolic powers and symbolic normal cone}   \label{subsec-valuation-symbolic-power}
It is natural here to mention the notion of symbolic powers of an ideal. Let $\p \subset R$ be a prime ideal. For $n \geq 0$, the {\it $n$-th symbolic power} $\p^{(n)}$ is by definition the smallest $\p$-primary ideal that contains $\p^n$, also for $n \leq 0$ we put $\p^{(n)} = R$ (see \cite[Section 3.9]{Eisenbud} for more about the notion of symbolic power). One sees that the symbolic powers $(\p^{(n)})_{n \in \Z}$ form a multiplicative filtration. We denote the corresponding Rees algebra and associated graded by $\mathcal{A}_{\p^{(*)}}(R)$ and $\gr_{\p^{(*)}}(R)$ and call them the {\it symbolic Rees algebra} and {\it symbolic associated graded of $\p$} respectively. In analogy with the usual normal cone, we call $\Spec(\gr_{\p^{(*)}}(R))$, the {\it symbolic normal cone} to the subscheme defined by $\p$ in $\Spec(R)$. Now, suppose that the local ring $R_\p$ is a discrete valuation ring and let $v$ be the corresponding valuation. One shows that the $n$-th symbolic power is given by the $n$-th subspace in the filtration $\F_v$, that is, $\p^{(n)} = F_{v, n} = \{ f \in R \mid v(f) \geq n\}$. %Let $X = \Spec(R)$ and $Y = V(\p)$ the subscheme defined by $\p$. 
%In analogy with the usual normal cone, we call $\Spec(\gr_v(R))$, the {\it symbolic normal cone} to the subscheme defined by $\p$ in $\Spec(R)$.

It is important to point out that the algebra $\mathcal{A}_{\p^{(*)}}(R)$ could be non-finitely generated. Some interesting examples of this situation can be found in \cite{Cutkosky}. Lemma \ref{lem-main} shows finite generation of $\mathcal{A}_{\p^{(*)}}(R)$ for certain prime ideals of height $1$. 

To illustrate the above concepts, below we give an example of the associated graded with respect to powers of a prime ideal versus the associated graded with respect to its symbolic powers. To facilitate the computation we state a lemma first. Let $\F$ be a filtration on $R$. Take $f \in R$ and let $n = v_\F(f)$, that is, $f \in F_n$ but $f \notin F_{n+1}$. Provided that $n \neq \infty$, we denote the image of $f$ in $F_n / F_{n+1} \subset \gr_\F(R)$ by $\bar{f}$ and call it the {\it initial form of $f$} with respect to the filtration $\F$. The next lemma is well-known and straightforward to prove (see \cite[Exercise 5.3]{Eisenbud}).

\begin{lemma} \label{lemma-initial-ideal} Let $R$ be an algebra and $J \subset I$ ideals in $R$. Then the associated graded $\gr_{I/J} (A/J)$ is naturally isomorphic to the quotient algebra $\gr_{I} (A)/\In(J)$, where $\In(J)$ is the ideal generated by the initial forms of elements of $J$ with respect to the filtration given by the powers of $I$.
\end{lemma}

\begin{example}[Normal cone and symbolic normal cone for an elliptic curve]\label{elliptic curve} Let $S = \k[x, y, z]$ and consider the algebra $R = S/(f)$ where $$f(x, y, z) = y^2 z - x^3 + xz^2.$$
The algebra $R$ is the homogeneous coordinate ring of the elliptic curve $X \hookrightarrow \mathbb{P}^2$ defined by the homogeneous polynomial $f$. Let $\tilde{x}, \tilde{y}, \tilde{z}$ be the images of $x, y, z$ in $R$.
Consider the prime ideal $\p = \sqrt{(\tilde{z})} = (\tilde{x}, \tilde{z}) \subset R$. Note that $\gr_{(x, z)}(S)$ is naturally isomorphic to $S$. We have $f \equiv 0 \text{ mod }(x, z)$ and $f \equiv y^2 z \text{ mod }(x, z)^2.$ This implies that the initial form $\bar{f}$ of $f$, with respect to the ideal $\p$, is $y^2 z$. By Lemma \ref{lemma-initial-ideal}, we have:
$$\gr_{\p}(R) \cong S/(y^2 z).$$
%Geometrically, as divisors on the affine space $\mathbb{A}^3$, we have:
%$$\Spec(\gr_{\p} R) = 2 \{ y = 0 \} + \{ z = 0 \}.$$
Let $\gr_{\p^{(*)}}(R) = \bigoplus_{n=0}^{\infty} \p^{(n)}/\p^{(n+1)}$ be the symbolic associated graded of $\p$. It is $(\Z_{\geq 0} \times \Z_{\geq 0})$-graded with the first $\Z_{\geq 0}$-grading inherited from $R$ and the second one given by the direct sum.
%By Lemma \ref{symbolic gr f-gen}, it is finitely generated: indeed, we can take $s = y^2, k = 0$ (TODO: more details.)
We have $\tilde{y}^2 \tilde{z} = \tilde{x}^3 - \tilde{x} \tilde {z}^2 \in \p^{3}$ and so $\tilde{z} \in \p^{(3)}$. One observes that $\tilde{z} \not\in \p^{(4)}$ (because $\tilde{z}$ vanishes of order $3$ at the point $(0:1:0)$ on $X$) and so with respect to the $(\Z_{\geq 0} \times \Z_{\geq 0})$-grading, the bidegree of $\tilde{z}$ is $(1, 3).$ Similarly, the bidegrees of $\tilde{x}$ and  $\tilde{y}$ are $(1, 1)$ and $(1, 0)$ respectively.
%Note: geometrically, the calculations correspond to the calculations of order-of-vanishing. For example, we have: $\tilde{z} = 0 \Rightarrow \tilde{x}^3 = 0$; i.e., $\tilde{z}$ vanishes to order 3 at $(0 : 1 : 0)$.
We now show that the images of $\tilde{z}, \tilde{x}, \tilde{y}$ generate the algebra $\gr_{\p^{(*)}}(R)$. Let $R' \subset \gr_{\p^{(*)}}(R)$ be the subalgebra generated by the images of $\tilde{x}, \tilde{y}, \tilde{z}$. On the one hand, for $n \geq 0$, we have:
\begin{align*}
\dim_k(R'_n) &= |\{ n_1 (1, 3) + n_2 (1, 1) + n_3 (1, 0) \mid n_1 + n_2 + n_3 = n, \, n_i \ge 0 \}| \\
&= 3n.
\end{align*} On the other hand, using the Riemann-Roch theorem, one computes that $\dim_\k (\gr_{\p^{(*)}} (R))_n = \dim_\k R_n = 3n$. This shows that $R' = \gr_{\p^{(*)}}(R)$ as claimed. Hence we have:
$$\gr_{\p^{(*)}}(R) \cong S/(y^2 z - x^3).$$
Geometrically, $X$ degenerates to the cuspidal cubic curve $y^2 z = x^3$ (cf. \cite[Example 4.2]{Anderson}).
\end{example}

\begin{remark}[Balanced normal cone] 
Although not needed here in this paper, we would like to mention another variant of the normal cone, namely the {\it balanced normal cone} (see \cite{Knutson}). Let $R$ be a Noetherian ring and $I \subset R$ an ideal. Let $v: R \to \Z_{\geq 0} \cup \{\infty\}$ be the quasivaluation associated to the filtration by powers of $I$, that is, for $f \in R$, $v(f) = n$ if $f \in I^n \setminus I^{n+1}$. If there is no such $n$, $v(f) = \infty$. Then define $\overline{v}$ by $$\overline{v}(f) = \lim_{n \to \infty} \frac{v(f^n)}{n}.$$
Samuel proved that the above limit exists and thus $\overline{v}(f)$ is well-defined. It is also known that $\overline{v}$ has values in 
$\Q$. One sees that $\overline{v}$ is a quasivaluation on $R$. Hence one can form the associated graded ring:
$$\overline{\gr}_I(R) = \bigoplus_{n \in \mathbb{Q}} \{ x \in R \mid \overline{v}(f) \geq n \}/\{ f \in R \mid \overline{v}(f) > n \}.$$
Nagata showed that $\sup_{f \in R} \mid \overline{v}(f) - v(f) | < \infty.$ From this and the inequality $\overline{v}(f) \geq v(f)$, $\forall f \in R$, it follows that $\overline{\gr}_I(R)$ is a Noetherian ring. The scheme $\Spec(\overline{\gr}_I(R))$ is called the {\it balanced normal cone} for the subscheme of $\Spec(R)$ defined by $I$. Two important properties of $\overline{\gr}_I(R)$ are the following: (1) It is a reduced ring. (2) The kernel of the natural homomorphism $\gr_I(R) \to \overline{\gr}_I(R)$
coincides with the nilradical of $\gr_I(R)$, and moreover, when the kernel vanishes, we have $\gr_I(R) = \overline{\gr}_I(R)$.
\end{remark}

\subsection{Khovanskii bases and flat degenerations} \label{subsec-Khovanskii-bases}
Let $R$ be a finitely generated positively graded domain and let $\v: R \setminus \{0\} \to \Z^r$ be a homogeneous valuation.

\begin{definition}[Khovanskii basis \cite{Kaveh-Manon-NOK}]     \label{def-Khov-basis}
A subset $\B \subset R$ is called a {\it Khovanskii basis} for $(R, \v)$ if the image of $\B$ in $\gr_\v(R)$ forms a set of algebra generators. 
\end{definition}

One can show that if $\v$ has full rank and the base field is algebraically closed, $\B$ is a Khovanskii basis if $\v(\B) = \{\v(f) \mid f \in \B\}$ generates $S(R, \v)$ as a semigroup.  

Given a Khovanskii basis one can do algebra operations in $R$ algorithmically, in particular one can represent any element of $R$ as a polynomial in the Khovanskii basis elements using a simple algorithm known as the {\it subduction algorithm}. In \cite{Kaveh-Manon-NOK} a theory of Khovanskii bases is developed. Moreover, it is shown that a finite set of algebra generators $\B$ is a Khovanskii basis, with respect to some valuation, if and only if the tropical variety of the ideal of relations among the elements of $\B$ contains a ``prime cone''. 
 
Whenever we have a Khovanskii basis, we can construct a deformation of $R$ to the associated graded $\gr_\v(R)$. 
%If moreover, $R$ is positively graded and the valuation $\v$ is homogeneous with respect to the grading, we get a degeneration of $\Proj(R)$ to the projective variety $\Proj(\gr_\v(R))$.
More precisely we have the following (\cite[Proposition 5.1]{Anderson} and \cite[Proposition 2.2]{Teissier}).

\begin{theorem}   \label{th-Anderson}
Suppose $(R, \v)$ has a finite Khovanskii basis, then 
there is a finitely generated, positively graded, flat $\k[t]$-subalgebra $\mathcal{A} \subset R[t]$, such that:
\begin{itemize}
\item[(a)] $\mathcal{A}[t^{-1}] \cong R[t, t^{-1}]$ as $\k[t, t^{-1}]$-algebras, and
\item[(b)] $\mathcal{A} / (t) \cong \gr_\v(R)$.
\end{itemize}
\end{theorem}

%If moreover we assume that $R$ is positively graded and the valuation $\v$ is homogeneous with respect to this grading, 
We have the following geometric interpretation of Theorem \ref{th-Anderson}.
\begin{corollary}  \label{cor-Anderson}
%With assumptions as above, 
There is a degeneration of $X = \Proj(R)$ to the variety $\Proj(\gr_\v(R))$ (in the sense of Definition \ref{def-intro-flat-degen}).
\end{corollary}
When $\v$ is full rank and the base field $\k$ is algebraically closed, the associated graded $\gr_\v(R)$ is isomorphic to the semigroup algebra $\k[S(R, \v)]$ and we get a toric degeneration of $\Proj(R)$ to the toric variety $\Proj(\k[S(R, \v)])$.

{Below we prove a converse to Theorem \ref{th-Anderson}. Roughly speaking, it says that %if a toric degeneration is such that the standard $\mathbb{G}_m$-action on the base $\A^1$ lifts to the total space of the family then the 
any toric degeneration must come from a full rank valuation. Its proof relies on \cite[Theorem 4]{Kaveh-Manon-NOK}.}
\begin{theorem}   \label{th-degen-valuation}
Let $R$ be a positively graded domain and let $\mathcal{A}$ be a finitely generated positively graded $\k[t]$-module and domain with the following properties:
\begin{itemize}
\item[(a)] $\mathcal{A}[t^{-1}] \cong R[t, t^{-1}]$ as $\k[t]$-modules and graded algebras.
\item[(b)] The algebra $R' = \mathcal{A}/(t)$ is a graded semigroup algebra $\k[S]$ where $S \subset \Z_{\geq 0} \times \Z^d$ is a finitely generated semigroup. 
\item[(c)] The standard $\mathbb{G}_m$-action on $\k[t]$ extends to an action on $\mathcal{A}$ respecting its grading. Moreover, this $\k[t]$-action acts through $\mathbb{G}_m^d$ acting on the semigroup algebra $R'$.
\end{itemize}
Then there is a full rank valuation $\mathfrak{v}: R \setminus \{0\} \to \Z_{\geq 0} \times \Z^d$ such that $R' \cong \k[S(R, \mathfrak{v})]$. 
\end{theorem}

\begin{proof}
{Let $\mathcal{A} = \bigoplus_{n \in \Z} F_i$ be the isotypical decomposition of $\mathcal{A}$ with respect to the $\mathbb{G}_m$-action.  The parameter $t$ is a non-zero divisor, so we must have $tF_i \subset F_{i+1}$.   By setting $t = 1$ we obtain the isomorphism of algebras $R \cong \mathcal{A}/(t-1)$, furthermore both of these algebras are identified with $\cup_{n \in \Z} F_n$ as vector spaces.  In this way, each $F_n$ is mapped isomorphically to a subspace of $R$, and these spaces form a homogeneous filtration $\mathcal{F}$ with $R' = \gr_{\mathcal{F}}(R)$. 

We select generators $u_1, \ldots, u_k \in R'$ which are homogeneous with respect to the $S$-grading on $R'$; this in turn implies they are homogeneous with respect to both the positive grading on $R'$ inherited from $\mathcal{A}$ and the action of $\mathbb{G}_m$.  If $u_i \in F_j/F_{j+1} \subset R'$ we select a homogeneous lift $y_i \in F_j$; we claim the $y_j$ and $t$ generate $\mathcal{A}$.  Let $F_{n, d} \subset F_n$ be the homogeneous degree $d$ part of $F_n$.  We must have $F_n = \bigoplus_d F_{n, d}$, and the $F_{n, d}$ define a vector space filtration of the $d$-graded component of $R$. It suffices to show that every $f \in F_{n, d}$ can be written as a polynomial in the $y_j$. For this we use a variant of the subduction algorithm (see \cite[Algorithm 2.11]{Kaveh-Manon-NOK}).  Since the $u_j$ generate $R' = \bigoplus_n F_n/F_{<n}$, we can find a monomial term $m = C_{\alpha}y^{\alpha}$ such that $f - m \in tF_{n-1, d}$.  Now we can repeat this procedure with $f-m$, noting that it must terminate eventually as $A$ is positively graded. 

It now follows that the images $z_j \in R$ of the $y_j$ generate $R \cong \mathcal{A}/(t-1)$. We now consider a monomial weighting $\r = (r_1, \ldots, r_k)$ of these generators, where $z_j$ is given the weight $r_j$ such that $u_j \in F_{r_j}/F_{r_j -1}$. Let $J$ be the homogeneous ideal which vanishes on the $z_j$. We consider $\r$ as a point in the Gr\"obner fan of $J$. Choose a Gr\"obner basis $G \subset J$ accordingly, and let $J_\r \subset k[x_1, \ldots, x_k, t]$ be the ideal which cuts out the Gr\"obner degeneration of $J$ corresponding to $\r$.   Observe that each polynomial in $J_\r$ must vanish on the $y_k$ and $t$ in $\mathcal{A}$.  But $k[x_1, \ldots, x_k, t]/J_\r$ has the same Krull dimension as $\mathcal{A}$, and both algebras are domains, so it must follow that $\mathcal{A} = k[x_1, \ldots, x_k, t]/J_\r$.  

Let $I$ be the binomial ideal which vanishes on the $u_j \in R'$. We have $\In_\r(J) = I$, and that $\r$ is in the lineality space of $I$. Now a standard argument from Gr\"obner theory shows that if $B$ is a sufficiently small relatively open ball in the lineality space of $I$, we must have $\In_{\r + s}(J) = \In_\r(J) = I$, for any $s \in B$.  Since $I$ is prime, this implies that the Gr\"obner fan of $J$ contains a prime cone $C$ of dimension $d$.  Now \cite[Theorem 4]{Kaveh-Manon-NOK} implies that there is a full rank valuation $\mathfrak{v}$ on $R$ with $R' = \k[S(R, \mathfrak{v})]$.}
\end{proof}

\section{Degeneration in stages}  \label{sec-degen-stages}
In this section we prove the following which is one of the main results of the paper.
\begin{theorem}   \label{th-degen-in-stages}
Let $R$ be a finitely generated positively graded $\k$-domain. Let $X = \Proj(R)$ and $d = \dim(X)$. Then there exists a sequence of $\k$-domains $R_0 =R, \ldots, R_{d-1}$ and degenerations $X_i = \Proj(R_i) \rightsquigarrow X_{i+1}=\Proj(R_{i+1})$, $i=0, \ldots, d-2$, such that $X_{d-1}$ is a complexity-one $T$-variety for an action of torus $T=\G_m^{d-1}$.  
\end{theorem}

The theorem is a corollary of the following lemma. Its proof in turn relies on Lemma \ref{lem-main} proved in the appendix, using some  technical commutative algebra.  
\begin{lemma}   \label{lem-degen-in-stages}
Let $A$ be a finitely generated $(\Z_{\geq 0} \times \Z^r)$-graded $\k$-domain and let $X = \Proj(A)$ (where $\Proj$ is taken with respect to the $\Z_{\geq 0}$-grading). The $\Z^r$-grading on $A$ induces an action of the torus $T = \G_m^r$ on $X$. We assume that $d = \dim(X) \geq 2$, $0 \leq r \leq d - 2$ and the following hold:
\begin{itemize}
\item[(1)] $T$-stabilizer of a general point in $X$ is finite.
\item[(2)] $\dim(X^\us) \leq r$ where $X^\us$ denotes the unstable locus for the $T$-action on $X$.
\item[(3)] $\dim(A^T) \geq d+1 - r$.
\end{itemize}
Then $X$ can be degenerated to $X' = \Proj(A')$ (in the sense of Definition \ref{def-intro-flat-degen}) where $A'$ is a finitely generated $(\Z_{\geq 0} \times \Z^{r+1})$-graded $\k$-domain with the following properties: Consider the $T' = \G_m^{r+1}$-action on $X'$ induced from the $\Z^{r+1}$-grading, then (1)-(3) above hold for $(A', X', T', r+1)$ in place of $(A, X, T, r)$ respectively. 
\end{lemma}
\begin{proof}
We construct $A'$ as the symbolic associated graded of a height one prime ideal in $A$. Let $Y = X /\!/ T = \Proj(A^T)$ be the GIT quotient of $X$ with $\pi: X^{\ss} \to Y$ the GIT quotient map. 
Let $k > 0$ be sufficiently large so that the Veronese subalgebra $(A^T)^{[k]} = \bigoplus_{i \geq 0} A^T_{ki}$ is generated in degree $1$. Take a finite set $\{f_0, \ldots, f_s\}$ of homogeneous degree $1$ algebra generators for $(A^T)^{[k]}$ to get an embedding $Y \hookrightarrow \P^s$ and consider $\pi: X^\ss \to Y \to \P^s$. For $0 \neq f \in (A^T)_k$, let $H_f \subset Y$ be the hyperplane section defined by $f$. Since $\dim(A^T) \geq d+1-r$ we have $\dim(Y) \geq d-r \geq 2$. Thus we can apply the Bertini irreducibility theorem (see \cite[Theorem 0.1]{Bertini}) to conclude that if $f$ is in general position then $\pi^{-1}(H_f) \subset X^\ss$ is an irreducible subvariety. Moreover, $\pi^{-1}(H_f)$ intersects the open subset $U \subset X^\ss$ of points with finite $T$-stabilizer. 
By assumption $\dim(X \setminus X^\ss) = \dim(X^\us) \leq r < d - 1$. That is, $\textup{codim}(X^\us) \geq 2$. It follows that the subvariety in $X$ defined by the principal ideal $(f)$ is irreducible and hence coincides with $Z = \overline{\pi^{-1}(H_f)}$. We note that $\pi^{-1}(H_f)$ and hence $Z$ are $T$-invariant subvarieties. This then implies that the radical ideal $\p = \sqrt{(f)}$, which is the ideal of the subvariety $Z$, is a homogeneous $T$-invariant prime ideal in $A$. Moreover, since $f \in (A^T)_k$ is in general position, we can assume that:
\begin{itemize}
\item The local ring $A_\p$ is a discrete valuation ring.
\item $Z \cap U \neq \emptyset$
\end{itemize} 
Now let $A' = \gr_{\p^{(*)}} A$ be the symbolic associated graded of $\p$. It is a $(\Z_{\geq 0} \times \Z^r \times \Z)$-graded algebra where the last $\Z$-grading is the natural $\Z_{\geq 0}$-grading on $A' = \bigoplus_{i \geq 0} \p^{(i)} / \p^{(i+1)}$. Let $X' = \Proj(A')$. Let $T' = \G_m^{r+1}$, then the $\Z^{r+1}$-grading on $A'$ induces a $T'$-action on $X'$. 
Also the torus $T' = T \times \G_m$ acts on $Z \times \A^1$ where $T$ acts on $Z$ and $\G_m$ acts on $\A^1$ in the standard way. 

By Lemma \ref{lem-main} we know that there is a natural embedding of the polynomial ring $(A/\p)[u]$ into $A'$ and $A'$ is a finite module over $(A/\p)[u]$. Moreover, this embedding induces a finite morphism $\phi: \Proj(A') \to \Proj((A/\p)[u])$. One verifies that the embedding $(A/\p)[u] \hookrightarrow A'$ preserves the $(\Z_{\geq 0} \times \Z^{r+1})$-gradings and hence $\phi$ is $T'$-equivariant. 

%We note that $A'^{T'} = (A/\p)^T$ and hence $(0,0) \in \Z^r \times \Z$ is not a regular value for the $T'$-action on $(X', A')$ defined above. In other words, the GIT quotient $X'  /\!/_{(0, 0)} T'$ does not have correct dimension. On the other hand, $(0, 1)$ is a regular value and we can consider the GIT quotient $X'  /\!/_{(0,1)} T'$.
%This amounts to twisting the $T'$-action by the character $(0, 1)$. Twisting the $T'$-action on $A'$ (respectively $(A/\p)[u]$) does not change the action of $T'$ on $X'$ (respectively $\Proj((A/\p)[u])$).

The variety $\Proj((A/\p)[u])$ is the projectivization of $\tilde{Z} \times \A^1$, where $\tilde{Z}$ is the affine cone over $Z$. It is the product $Z \times \A^1$ with a point $\infty$ added. The point $\infty$ is in the closure of each fiber $\A^1$. It is easy to see that the dimension of unstable locus for the $T'$-action is at most dimension of $T$-unstable locus on $Z$ plus $1$ which in turn is less than or equal to $r + 1$. 

Since $Z \cap U \neq \emptyset$ we know that the generic $T'$-stabilizer in $Z \times \A^1$ is finite. From $T'$-equivariance of $\phi$ it follows that the generic $T'$-stabilizer in $X'$ is also finite. %Let $Z^\us$ and $X^\us$ be the unstable locus for $Z$ for the action of $T$. Since $Z^\us \subset X^\us$ we know that $\dim(Z^\us) \leq \dim(X^\us) < \dim(X) - 1$. On the other hand, one verifies that the unstable locus of $T'$-action on $Z \times \A^1$ is also $Z^\us$ and thus has codimension less than $1$. 
%We consider the twist the $T'$ by a character. It does not change the action on $Z'$ and $X'$.

Now by \cite[Chap I, Section 5]{Mumford} we know that inverse image of the $T'$-unstable locus of $\Proj((A/\p)[u])$ under the finite morphism $\phi$ contains the $T'$-unstable locus of $X'$. It follows that the $T'$-unstable locus of $X'$ has dimension $\leq r+1$. To finish the proof we need to verify that $\dim(A'^{T'}) \geq d+1 - (r+1)$. We note that $A'^{T'} = A^T/ \p^T$ with $\p^T = \p \cap A^T$. Since $\p = \sqrt{fA}$ and $f$ is $T$-invariant, we have $\p \cap A^T = \sqrt{f A^T}$. Thus, $\p^T = \p \cap A^T$ is a height $1$ prime in $A^T$ and $\dim(A^T / \p^T) = \dim(A^T) - 1 \geq d+1 - r - 1$ as required. 
\end{proof}

\begin{proof}[Proof of Theorem \ref{th-degen-in-stages}]
Starting from $r=0$ repeatedly apply Lemma \ref{lem-degen-in-stages} until $r= \dim(X) - 2$. We arrive at a sequence of degenerations of $X$. 
\end{proof}

\begin{remark}
Observe that the Bertini theorem breaks down in dimension $1$, so our method above cannot produce a toric degeneration in general.
\end{remark}

\section{Curve case}    \label{sec-curve-case}
In this section we consider the case where the variety $X$ is a projective curve. We give an example of the homogeneous coordinate ring of a smooth projective curve such that for any choice of a homogeneous full rank valuation the corresponding associated graded (equivalently the value semigroup) is non-finitely generated. Also, as a corollary of Lemma \ref{lem-main} we see that given a (not necessarily smooth) projective curve $X$, we can find a very ample line bundle $L$ and a valuation $\v$ (corresponding to a smooth point on $X$) such that the associated graded of the homogeneous coordinate ring of $(X, L)$ with respect to $\v$ is finitely generated and hence $(X, L)$ has a toric degeneration. 

%The deformations of curves are very well studied. %in particular in the context of test configurations. 
%We can mention the book \cite{Pinkham} in this regard.

For a smooth point $p \in X$ let us consider the homogeneous valuation $\v_p: R \setminus \{0\} \to \Z_{\geq 0} \times \Z$ as follows. For every $f \in R$ put:
$$\v_p(f) = (\deg(f), \ord_p(f)).$$
(The first factor $\Z_{\geq 0}$ is equipped with the reverse ordering of integers). 

\begin{prop}[Finitely generated semigroup case] 
\label{prop-curve-finitely-gen}
Let $X$ be a (not necessarily smooth) projective curve. Let $L$ be a very ample line bundle on $X$ such that the divisor class of $L$ is a multiple of a smooth point $p \in X$. Then the algebra of sections $R=R(X, L)= \bigoplus_{i} H^0(X, L^{\otimes i})$ has a finite Khovanskii basis with respect to the valuation $\v_p$, i.e. the value semigroup $S(R, \v_p)$ is finitely generated. 
\end{prop}
\begin{proof}
From Lemma \ref{lem-main} it follows that $\gr_{\v_p}R$ is finitely generated which is what the proposition claims.
\end{proof}
A variant of Proposition \ref{prop-curve-finitely-gen} is observed in \cite[Remark 3.17]{Ilten-Wrobel}, and an example of the situation in Proposition \ref{prop-curve-finitely-gen} is Example \ref{elliptic curve}.  Also see the discussion after Corollary 17 in \cite{AKL}. 
  
\begin{prop}[Non-finitely generated semigroup case]  
\label{prop-curve-not-finitely-gen}
Suppose the base field $\k$ is uncountable. Let $X$ be a smooth projective curve with genus $>1$ and let $L$ be a very ample line bundle on $X$ such that no tensor power of $L$ has a divisor which is a multiple of a point. Then for any choice of $p \in X$ the value semigroup $S(R, \v_p)$ is not finitely generated and hence $(R, \v_p)$ does not have a finite Khovanskii basis (and hence no toric degeneration for any choice of a homogeneous valuation).
\end{prop}
\begin{proof}
Take a point $p \in X$ with the corresponding valuation $\v_p$. We know that the cone of the value semigroup $S(R, \v_p)$ is the cone $C = \{(k, x) \mid 0 \leq x \leq k\delta \}$ where $\delta = \deg(L)$ is the degree of the line bundle $L$. By contradiction suppose that the value semigroup $S(R, \v_p)$ is indeed finitely generated. Then for some $m$ there should exists $f \in R_m$ such that $\ord_p(f) = m\delta$. This means that the divisor of the section $f$ is $m\delta p$ which contradicts the assumption, namely the line bundle $L^{\otimes m}$ does not have a divisor which is a multiple of a single point. 
\end{proof}

In fact, it is also possible to find rational curves with a projective coordinate ring which cannot have a finite Khovanskii basis. In \cite[Corollary 3.14]{Ilten-Wrobel} this is shown for the projective coordinate ring of a very general integral rational plane curve of degree $d > 3$. 

\section{Cohen-Macaulay case}   \label{sec-CM-case}
%{\bf Suppose $A$ is Cohen-Macaulay. Then construction of degeneration follows from a theorem of Rees.}
Throughout this section $R$ is the homogeneous coordinate ring of a $d$-dimensional projective variety $X$.
When the ring $R$ is Cohen-Macaulay (in other words, $X$ is arithmetically Cohen-Macaulay or ACM for short), there is a simple construction of a degeneration of $X$ to a complexity-one $T$-variety.

%Let $R$ be a commutative ring. 
Recall that a sequence $f_1, \ldots, f_r$, $f_i \in R$, is called a {\it regular sequence} if $f_1$ is not a zero divisor and for every $i=2, \ldots, r$, the image of $f_i$ in $R / (f_1, \ldots, f_{i-1})$ is not a zero divisor. If $R$ is an algebra over a field $\k$, it follows from the definition of a regular sequence that $f_1, \ldots, f_r$ are linearly independent over $\k$.

Let $\p = (f_1, \ldots, f_r)$ be the ideal generated by the $f_i$. The following theorem is due to Rees. The idea of proof goes back to Macaulay and is based on a double induction (see \cite[Theorem 2.1]{Rees}). 
\begin{theorem}   \label{th-Rees}
Let $\bar{f}_i$ be the image of $f_i$ in $\p / \p^2$. Then the map $(R/\p)[t_1, \ldots, t_r] \to \gr_\p(R)$ which is identity on $R/\p$ and sends $t_i$ to $\bar{f}_i$, $i=1, \ldots, r$, gives an isomorphism between $(R/\p)[t_1, \ldots, t_r]$ and $\gr_\p(R)$.
\end{theorem}

%\begin{prop}  \label{prop-Rees-alg-fg}
%{\bf Consider the Rees algebra $\mathcal{A} = \mathcal{A}_\p(R) = \bigoplus_{i \geq 0} \p^i$. Then $\mathcal{A}$ is generated as an algebra by $\{1, f_1, \ldots, f_r\}$ where $1 \in \k=\mathcal{A}_0$ and $f_i \in \p = \mathcal{A}_1$.}
%\end{prop}
%\begin{proof}
%{\bf From definition.}
%\end{proof}

%\begin{prop}
%The associated graded $\gr_v(R)$ is isomorphic to the polynomial algebra $(R/\p)[t_1, \ldots, t_r]$. The isomorphism is given by sending $R/\p$ to $R/\p$ and $t_i$ to the image $\bar{f}_i$ of $f_i$ in $\p/\p^2$.
%It follows that $\gr_\v(R)$ is isomorphic to $\gr_\p(R)$.
%\end{prop}

Let $X$ be a projective variety of dimension $d \geq 2$. Fix an embedding of $X$ in a projective space $\P^N$ and let $R = \k[X]$ be its homogeneous coordinate ring. We write $R = \bigoplus_{i \geq 0} R_i$ where $R_i$ is the $i$-th homogeneous piece of $R$ with respect to its natural grading. We also let $\tilde{X} \subset \A^{N+1}$ denote the affine cone over $X$. %We let $z_0, \ldots, z_N$ be the coordinates on $\A^{N+1}$. Regarded as elements of $R$, the $z_i$ have degree $1$ and generate $R$ as an algebra.

%The following is a corollary of the Bertini irreducibility theorem. Cite Benoit's paper?. Add Takuya's text here and add a reference to the book Local cohomology in 24 hours.

%We recall that a projective subvariety is irreducible if and only if its corresponding homogeneous ideal is prime. 
The following is a corollary of the Bertini irreducibility theorem and unmixedness of Cohen-Macaulay rings. 
\begin{prop}   \label{prop-reg-seq-exists}
Suppose $R$ is Cohen-Macaulay. Then there exists a Zariski open subset $\mathcal{U} \subset R_1 \times \cdots \times R_1$ ($d-1$ times) such that for every $(f_1, \ldots, f_{d-1}) \in \mathcal{U}$ and any subset $\{i_1, \ldots, i_k\} \subset \{1, \ldots, d-1\}$ the ideal generated by $f_{i_1}, \ldots, f_{i_k}$ is prime.  
In particular, $(f_1, \ldots, f_{d-1})$ is a regular sequence in $R$ and  the ideal $\p = (f_1, \ldots, f_{d-1})$ is prime.  
\end{prop}
%\begin{proof}
%\end{proof}

Now fix a regular sequence $f_1, \ldots, f_{r} \in R_1$, $r \leq d$, such that $\p=(f_1, \ldots, f_{r})$ is prime.
Let $Y \subset X$ be the subvariety defined by the homogeneous prime ideal $\p$ in $X$. Also let $\tilde{Y} \subset \tilde{X}$ be the affine cone over $Y$ and let $\tilde{X}_0 = C_{\tilde{Y}}\tilde{X}$ denote the normal cone of $\tilde{Y}$ in $\tilde{X}$. Let $\pi: \tilde{\X} \to \A^1$ be the deformation to normal cone family which deforms $\tilde{X}$ to $\tilde{X}_0$ and let $\X$ be the projectivization of $\tilde{\X}$, that is, $\X = \Proj(\bigoplus_{i \geq 0} \p^i)$. The normal cone $\tilde{X}_0$ is the $\Spec$ of the associated graded algebra $\gr_\p(R)$. %One verifies that since $\p$ is generated by homogeneous elements of degree $1$, 
Since each graded piece of $\gr_\p(R)$ is a finite dimensional $\k$-vector space, $X_0 = \Proj(\gr_\p(R))$ is a projective variety.
%The ideal $\p$ is homogeneous and thus the graded algebra $\gr_\p(R)$ inherits a $\Z_{\geq 0}$-grading from that of $R$. For $m \geq 0$, the $m$-th graded piece $(\gr_\p(R))_m$ is:
%$$(\gr_\p(R))_m = \bigoplus_{i \geq 0} (\p^i)_m / (\p^{i+1})_m,$$
%where $(\p^i)_m = \p^i \cap R_m$, for all $i$. Note that for each $m \geq 0$, $\dim_\k(R_m) < \infty$. Since $\p$ is generated by elements of degree $1$ it follows that $(\p^i)_m = \{0\}$ for $i > m$. In particular, $\dim_\k((\gr_p(R))_m)  < \infty$ for all $m \geq 0$. This implies that $X_0 = \Proj(\gr_\p(R))$ is a projective variety. 

Let us give a precise description of the projective variety $X_0$. Since $\{f_1, \ldots, f_r\} \subset R_1$ is linearly independent over $\k$ we can extend it to a basis $\{f_1, \ldots, f_{N}\}$ for $R_1$. We observe that $\gr_\p(R) = \bigoplus_{i \geq 0} \p^i / \p^{i+1}$ is generated as an algebra by $\{1, \bar{f}_1, \ldots, \bar{f}_{N}\}$ where $\bar{f}_i$ is the image of $f_i$ in $\gr_\p(R)$. This choice of generators gives rise to an embedding of $X_0 = \Proj(\gr_\p(R))$ into $\P^N$. 
%the family $\X$ in $\P^{N} \times \A^1$ and we have a commutative diagram:
%$$\xymatrix{
%\X ~ \ar[rd]_{\pi} \ar@{^{(}->}[r] & \P^{N} \times \A^1 \ar[d]^{\text{pr}} \\ & \A^1\\
%}$$
%where $\text{pr}$ denotes the projection on the second factor. 
%Recall that $z_0, \ldots, z_N$ are degree $1$ algebra generators for $R$ corresponding to the coordinate functions in $\A^{N+1}$. Let $x_i$ be the image of $z_i$ in $R/\p$. Then $\{x_0, \ldots, x_N, t_1, \ldots, t_r\}$ is an algebra generating set for $(R/\p)[t_1, \ldots, t_r]$ and in view of Theorem \ref{th-Rees}, $\{x_0, \ldots, x_N, \bar{f}_1, \ldots, \bar{f}_r\}$) is an algebra generating set for $\gr_\p(R)$. 
%This homogeneous algebra generating set for $\gr_\p(R)$ gives an embedding $X_0 = \Proj(\gr_\p(R)) \hookrightarrow \P^{N+r}$. 
Consider the rational map $p: \P^{N} \dashrightarrow \P^{N-r}$ given by: 
$$(x_1: \cdots x_{N+1}) \mapsto (x_{r+1}: \cdots : x_{N+1}).$$
This rational map is regular outside the locus $Z$ where $x_{r+1} = \cdots = x_{N+1} = 0$. Also $p: \P^{N} \setminus Z \to \P^{N-r}$ is a vector bundle and fiber over each point can be identified with $\A^r$. 
 Moreover, $p: X_0 \setminus Z \to Y$ is also a vector bundle with fibers isomorphic to $\A^r$ (Theorem \ref{th-Rees}). We thus obtain the following: 
\begin{prop}   \label{prop-desc-X_0}
The variety $X_0 = \Proj(\gr_\p(R))$ is the disjoint union of a vector bundle of rank $r$ over $Y$ and a copy $P$ of the projective space $\P^{r}$. Moreover, each fiber of the vector bundle union with $P \cong \P^{r}$ is isomorphic to $\P^{r+1}$. Also $X_0$ has a natural action of torus $T=\G_m^{r}$. It acts on the fibers of the vector bundle and on the projective space $P$ in the natural way. Thus when $r=d-1$, $X_0$ is a complexity-one $T$-variety.
\end{prop}

\begin{remark}   \label{rem-CM-case-Khov-basis}
With notation as above, let $\v: R \setminus \{0\} \to \Z_{\geq 0} \times \Z^{d-1}$ be the homogeneous valuation corresponding to the regular sequence $(f_1, \ldots, f_{d-1})$. We obtain a Khovanskii basis for $(R, \v)$ by augmenting any generating set of $R$ with the regular sequence $f_1, \ldots, f_{d-1}$.
\end{remark}

%\begin{remark}   \label{rem-X-0-not-smooth}
%Unfortunately the vairety $X_0$ is almost never smooth as can be seen by looking at the affine charts around $P \cong \P^{r}$. For $1 \leq j \leq r$ consider the affine open subset $\tilde{U}_j \subset X_0$ defined by $t_j \neq 0$. We have:
%\begin{equation} \label{equ-tilde-U-j}
%\tilde{U}_j = \Spec(\bigcup_{m \geq 0} \{ \frac{p}{t_j^m} \mid p \in (R/\p)[t_1, \ldots, t_r] \text{ homogeneous of degree } m\}).
%\end{equation} 
%It is straightforward to verify that the righthand side of \eqref{equ-tilde-U-j} is isomorphic to $(R/\p)[t_1, \ldots, \hat{t}_j, \ldots, t_r]$ (where $\hat{t}_j$ indicates that $t_j$ has been removed). It follows that $\tilde{U}_j \cong \tilde{Y} \times \A^{r-1}$. Since $\tilde{Y}$ is the affine cone over the projective curve $Y$ it is almost never smooth.
%\end{remark}

Now we assume that $R$ is both Cohen-Macaulay and normal. %in particular we assume $R$ is generated by its degree $1$ component.  
A general principle in algebraic geometry dictates that the complete intersection $Y = \Proj(R/\p)$ inherits nice properties from $X$, so $Y$ is likewise projectively normal and ACM (see e.g. \cite[Corollary 3.5.6]{FOV}).   This implies that $Y$ is a smooth curve equipped with an embedding $i:Y \hookrightarrow \P^N$.  Let $\mathcal{L} = i^*\O(1)$ be the corresponding very ample line bundle on $Y$; we have $R/\p \cong \bigoplus_{n \geq 0} H^0(Y, \mathcal{L}^{\otimes n})$.   With these assumptions, Proposition \ref{prop-curve-finitely-gen} and Theorem \ref{th-Rees} imply the following. 

\begin{prop}\label{prop-CM-norm-toric}
If for some $n > 0$ the divisor class of $\mathcal{L}^{\otimes n}$ is a multiple of a point $p \in Y$, then $R$ can be degenerated to a finitely generated affine semigroup algebra (in the sense of Theorem \ref{th-Anderson}). 
\end{prop} 

Let $h_R(n)$ be the Hilbert polynomial of $R$.  We write $h_R(n)$ as a sum of binomials:

\begin{equation}
h_R(n) = \sum_{i = 0}^d c_i \binom{ n + i}{i},  \ \ \ c_0, \ldots, c_d \in \Z.\\
\end{equation}

\noindent
The graded ring $R/\p$ is obtained from $R$ by a complete intersection of $d-1$ linear forms, as a consequence we have:

\begin{equation}\label{eq-curve-hilbert}
h_{R/\p}(n) = c_d(n + 1) + c_{d-1}, \\
\end{equation}

\noindent
and the genus of $Y$ is $ 1- c_d - c_{d-1}$. 

\begin{corollary}\label{cor-hilbert-degen}
Let $R$ and $h_R(n)$ be as above, then if $c_d + c_{d-1}$ equals $0$ or $1$, the algebra $R$ can be degenerated to a finitely generated affine semigroup algebra.  In particular, a normal hypersurface of degree $2$ or $3$ has a toric degeneration.  
\end{corollary}

\begin{proof}
If the genus of $Y$ is $0$ or $1$ then the divisor class of $\mathcal{L}$ is a multiple of a point, now we use Proposition \ref{prop-CM-norm-toric}. 
\end{proof}

%Perhaps a remark on a generic hypersurfaces of degree 2 and 3 here; these are normal, but the Newton-Polytope of such a polynomial is a regular simplex of sidelength 2 or 3.  The edges of these simplices are evidently not the newton polytopes of an irreducible polynomial.  This shows that a change of coordinates is necessary to find the toric degeneration. 

If $c_d + c_{d-1} = 1$ then $Y \cong \P^1$ and $\mathcal{L} \cong \O(\ell)$ for some $\ell \in \Z_{>0}$.  It is easy to see from \eqref{eq-curve-hilbert} that $\ell = c_d$, the degree of $X \subset \P^N$.  It follows that $R/\p$ is the $c_d$-th Veronese subring of the polynomial ring $\k[x, y]$.  Let $\Delta(\vec{1}, \ell)$ be the simplex with vertices $(0, \ldots, 0)$, $(1, \ldots, 0), \ldots, (0, \ldots, 1, 0)$, $(0, \ldots, 0, \ell)$, and let $S(\vec{1}, \ell)$ be the graded affine semigroup algebra associated to $\Delta(\vec{1}, \ell)$. From Theorem \ref{th-Rees} we can conclude the following. 

\begin{corollary}
Let $R$ be a normal Cohen-Macaulay domain as above, then the following are equivalent:\\

\begin{enumerate}
\item $R$ has a flat degeneration to $S(\vec{1}, c_d)$ (in the sense of Theorem \ref{th-Anderson}).
\item A complete intersection in $X$ by $d-1$ generic hyperplanes is isomorphic to $\P^1$.
\item $c_d + c_{d-1} = 1$.
\end{enumerate} 

\end{corollary}

The algebra $S(\vec{1}, c_d)$ is a quotient of the polynomial ring $\k[x_0, \ldots, x_{c_d}, y_1, \ldots, y_{d-1}]$ by the $2\times 2$ minors of the matrix:

$$
\begin{bmatrix}
x_0 & x_1 & \cdots & x_{c_d-1}\\
x_1 & x_2 & \cdots & x_{c_d}\\
\end{bmatrix}
$$
\smallskip

\noindent
These forms are a quadratic square-free Gr\"obner basis.  It follows that $R$ can likewise be presented as a quotient of $\k[x_0, \ldots, x_{c_d}, y_1, \ldots, y_d]$ by a quadratic square-free Gr\"obner basis.  The algebra $S(\vec{1}, c_d)$ is Gorenstein if and only if $c_d = 2$, this is also the case when $R$ can be presented as coordinate ring of a hypersurface. 
 
\begin{remark}\label{rem-ACM-toric}
Let $R$ be Cohen-Macaulay and normal, with $X$ and $X_0$ as above.  We choose a point $p \in Y$ and consider the ample divisor $D_p \subset X_0$ defined by the copy of $\mathbb{P}^{r+1}$ obtained from the fiber over $p$ as in Proposition \ref{prop-desc-X_0}.  The projective coordinate ring of $D_p$ is a polynomial ring over the projective coordinate ring of the divisor defined by $p \in Y$.  By Proposition \ref{prop-CM-norm-toric}, $X_0$ can be degenerated to a toric variety using this alternative embedding.  It follows that any normal projective variety with an arithmetically Cohen-Macaulay embedding can be degenerated in multiple steps to a toric variety.    
\end{remark}

\section{Existence of a finite Khovanskii basis and degeneration in one step}    \label{sec-finite-Khovanskii-basis}
In this section, using results from Section \ref{sec-degen-stages}, we show that a finitely generated positively graded algebra and domain with Krull dimension $d+1$ has a rank $d$ valuation with finitely generated associated graded algebra. This then implies that the projective variety corresponding to $R$ can be degenerated (in one step) to a complexity-one $T$-variety (in the sense of Definition \ref{def-intro-flat-degen}). 

More precisely, we have the following:
\begin{theorem}   \label{th-main}
Let $R$ be a finitely generated positively graded algebra and domain with $\dim(R) = d+1$ and $R_0 = \k$. Then there exists a homogeneous rank $d$ valuation $\mathfrak{v}: R \setminus\{0\} \to \Z_{\geq 0}  \times \Z^{d-1}$ such that the associated graded $\gr_\v(R)$ is finitely generated (in other words, $(R, \mathfrak{v})$ has a finite Khovanskii basis). 
\end{theorem}

\begin{remark}  \label{rem-fg-gr-vs-semigroup}
Note that the finite generation of $\gr_\v(R)$ as an algebra implies the finite generation of $S(R, \v)$ as a semigroup but not the other way around. Since the valuation $\v$ is not full rank, the finite generation of $\gr_\v(R)$ is stronger than the value semigroup be finitely generated.  
\end{remark}

\begin{corollary}   \label{cor-toric-degen-proj-var}
With notation as above, the projective variety $X = \Proj(R)$ has a degeneration to the complexity-one $T$-variety $\Proj(\gr_\v(R))$ (in the sense of Definition \ref{def-intro-flat-degen}). %Moreover, the standard action of $\G_m$ on $\A^1$ lifts to the total space of the degeneration.
\end{corollary}

%In fact, more generally, we can prove a multi-graded version of Theorem \ref{th-main}.
%\begin{theorem}
%Let $R$ be a finitely generated $\Z_{\geq 0}^r$-graded algebra and domain (for some $r > 0$) with $R_0 = \k$. Let $d$ be the transcendence degree of $0$-degree piece of field of fractions of $R$. Then there exists a  valuation $\mathfrak{v}: R \setminus\{0\} \to \Z^r_{\geq 0}  \times \Z^d$, where the $\Z_{\geq 0}^r$-component is just degree with respect to the grading of $R$, such that the value semigroup $S = S(R, \mathfrak{v}) \subset \Z_{\geq 0} \times \Z^d$ is finitely generated. 
%\end{theorem}

%\subsection{Construction of the valuation $\mathfrak{v}$} \label{subsec-construction-v}
The construction of the valuation $\mathfrak{v}$ relies on the following basic construction which we call {\it concatenation of valuations}. Let $R$ be as above and let $v: R \setminus \{0\} \to \Gamma_1$ be a valuation with values in an ordered group $\Gamma_1$. Recall that for $0 \neq f \in R$ with $v(f) = a$ we let $\bar{f}$ denote the image of $f$ in $R_{v \geq a} / R_{v > a} \subset \gr_{v}(R)$. The map $f \mapsto \bar{f}$ is a multiplicative homomorphism from $R \setminus \{0\}$ to the multiplicative set of homogeneous elements of $\gr_{v}(R)$. Also let $w: \gr_{v}(R) \setminus \{0\} \to \Gamma_2$ be a valuation on the associated graded of $v$ and with values in an ordered group $\Gamma_2$. We define a function $u: R \setminus \{0\} \to \Gamma_1 \times \Gamma_2$ as follows: for $0 \neq f \in R$ put:
\begin{equation}   \label{equ-concat-v-w}
u(f) = (v(f), w(\bar{f})).
\end{equation}

Let us equip $\Gamma_1 \times \Gamma_2$ with a lexicographic order. More precisely, for $(a_1, a_2)$, $(b_1, b_2) \in \Gamma_1 \times \Gamma_2$ we have $(a_1, a_2) > (b_1, b_2)$ if $a_1 > b_1$, or $a_1 = b_1$ and $a_2 > b_2$. The following is straightforward to verify.
\begin{prop}   \label{prop-concat-is-valuation}
With notation as above, the function $u: R \setminus \{0\} \to \Gamma_1 \times \Gamma_2$ is a valuation. 
\end{prop}
%\begin{proof}
%The function $u$ is well-defined, and clearly takes any nonzero element of the ground field $\k$ to $0$. To compute $u(fg)$, $0 \neq f, g \in R$, we must compute $(v(fg), w(\overline{fg})$.  But $\overline{fg} = \bar{f}\bar{g}$, so we have $(v(f) + v(g), w(\bar{f}) + w(\bar{g})) = u(f) + u(g)$.  To show the non-Archimedean property we consider the different possibilities for a sum $u(f + g)$.  If $v(f + g) > v(f) = v(g)$, then the lexicographic ordering on $\Gamma_1 \times \Gamma_2$ ensures that $u(f + g) > u(f) = u(g)$.  If $v(f + g) = v(f) < v(g)$ then $\overline{f + g} = \overline{f}$, so we obtain $u(f + g) = u(f) < u(g)$.  Finally, if $v(f + g) = v(f) = v(g)$ then $\overline{f + g} = \bar{f} + \bar{g}$, so in any case $u( f + g) \geq \min\{ u(f), u(g)\}$.  
%\end{proof}

We call the valuation $u$ the {\it concatenation of $v$ and $w$} and denote it by $u = v \circ w$.

The following proposition relates the associated graded algebra of the concatenation $v \circ w$ with those of $v$ and $w$.
\begin{prop}   \label{prop-gr-concat}
With notation as above, let $v:R\setminus\{0\} \to \Gamma_1$ be a valuation and $w: \gr_{v}(R)\setminus\{0\} \to \Gamma_2$ be a homogeneous valuation with respect to the $\Gamma_1$-grading on $\gr_{v}(R)$. Then we have an isomorphism of graded $(\Gamma_1 \times \Gamma_2)$-algebras:
\begin{equation}
\gr_{v \circ w}(R) \cong \gr_{w}(\gr_{v}(R)).
\end{equation}
\end{prop}
\begin{proof}
%We use the following notation:
%\begin{enumerate}
%\item $F$: the filtration on $R$ induced by $v$, with $F(\leq r) = v^{-1}(\leq r)$,
%\item $G$: the filtration on $gr_v(R)$ induced by $w$, with $G(\leq s) = w^{-1}(\leq s)$,
%\item $H$: the filtration on $R$ induced by $v \circ w$, with $H(\leq (r, s)) = v \circ w^{-1}(\leq (r, s))$.
%\end{enumerate}
Let $F = \F_{v}$, $G = \F_{w}$ be the filtrations on $R$ and $\gr_{v}(R)$ induced by $v$ and $w$ respectively. Also let $H = \F_{v \circ w}$ be the filtration on $R$ induced by the concatenation $v \circ w$. 
As $w$ is assumed to be homogeneous, the algebras $\gr_{v \circ w}(R)$ and $\gr_{w}(\gr_{v}(R))$ are both $(\Gamma_1 \times \Gamma_2)$-graded algebras, and $\gr_{w}(\gr_{v}(R))$ is the direct sum of the spaces $(G_{\geq s} \cap F_{\geq r}/F_{>r})/(G_{>s} \cap F_{\geq r}/F_{>r})$. We show that that the graded components of these algebras are isomorphic and that their homogeneous multiplication operations coincide. 

First we construct a map $\hat{\phi}: (G_{\geq s} \cap F_{\geq r}/F_{>r}) \to H_{\geq (r,s)}/H_{>(r, s)}$.  Any $\bar{f} \in F_{\geq r}/F_{>r}$ is represented by some $f \in F_{>r}$, we let $\hat{\phi}(\bar{f})$ for $\bar{f} \in (G_{\geq s} \cap F_{\geq r}/F_{>r})$ be the equivalence class of $f$ in $H_{\geq (r, s)}/H_{>(r, s)}$.  If we choose some other representative $f' \in F_{>r}$ then the difference $f - f'$ must lie in $F_{>r}$, and therefore $H_{>(r, s)}$, so this map is well-defined.  Any $g \in H_{\geq (r, s)}$ must have $v(g) \geq r$ and $w(\bar{g}) \geq s$ for $\bar{g} \in F_{\geq r}/F_{>r}$ by definition, it follows that  $\hat{\phi}$ is also onto. 

We claim that the kernel of $\hat{\phi}$ is $(G_{>s} \cap F_{\geq r}/F_{>r})$.  Clearly this space lies in the kernel of $\hat{\phi}$, and if $\hat{\phi}(\bar{f}) = 0$ it must be the case that a lift $f \in F_{\geq r}$ is in $H_{>(r, s)}$.  If $f \notin F_{>r}$ (i.e. $\bar{f} \neq 0$), then $w(\bar{f})$ will be strictly less than $s$.  It follows that the induced map on the quotient space $\phi: (G_{\geq s} \cap F_{\geq r}/F_{>r})/(G_{>s} \cap F_{\geq r}/F_{>r}) \to H_{\geq (r, s)}/H_{>(r, s)}$ is an isomorphism. 

Now let $\bar{f}_1, \bar{f}_2 \in \gr_{w}(\gr_{v}(R))$ be homogeneous. The product $\bar{f}_1\bar{f}_2$ is represented by the product $f_1f_2$, where $f_i$ is a lift of $\bar{f}_i$ in the definition of $\hat{\phi}$ above.  The construction above shows that multiplicative map $R \to \gr_{v \circ w}(R)$ is the composition of the maps $R \to \gr_{v}(R)$ and $\gr_{v}(R) \to \gr_{w}(gr_{v}(R)) \cong \gr_{v \circ w}(R)$, the latter being isomorphism as graded vector space, so it follows that $\phi(\bar{f}_1\bar{f}_2) = \phi(\bar{f}_1)\phi(\bar{f}_2)$.   
\end{proof}

\begin{remark}
For any valuation $w$ on $\gr_{v}(R)$, not necessarily homogeneous, there is an associated homogeneous valuation $\bar{w}$ defined as follows:
\begin{equation}
\bar{w}(\sum_\gamma f_{\gamma}) = \min\{w(f_{\gamma}) \mid f_{\gamma} \neq 0\}.
\end{equation}
Proposition \ref{prop-gr-concat} in fact shows that $\gr_{v \circ w}(R) \cong \gr_{\bar{w}}(\gr_{v}(R))$. 
\end{remark}

We might expect that $v \circ w$ has rank $\rank(v) + \rank(w)$, however the actual rank may be smaller, for example if $w$ is a grading function on $\gr_{v}(R)$ (see Section \ref{subsec-valuations}). When $w$ has rank $1$ this is all that can go wrong. 
\begin{prop}\label{concatenationrank}
Let $v: R \setminus \{0\} \to \Gamma$ have rank $r$, and let $w: \gr_{v}(R) \to \Q$ be homogeneous with respect to the $\Gamma$-grading on $\gr_{w}(R)$.  If $w$ is not a grading function then $v \circ w: R \setminus \{0\} \to \Gamma \times \Q$ has rank $r + 1$. 
\end{prop} 
\begin{proof}
By definition, the rank of $v \circ w$ is greater than or equal to $r$. Suppose by contradiction that it is equal to $r$. Then $S(R, v \circ w)$ lies in a hyperplane in $\Gamma \times \Q$ and projection onto the first factor in $\Gamma \times \Q$ gives an isomorphism from $S(R, v \circ w)$ to $S(R, v)$. Let $\pi: S(R, v) \to \Q$ be the linear function given by the projection onto the second factor in $\Gamma \times \Q$. We see that if $\bar{f} \in \gr_v(R)$ is homogeneous of degree $a$ then $w(\bar{f}) = \pi(a)$. But this implies that $w$ is a grading function. 
\end{proof}

%\begin{theorem}\label{concatenatebyrank1}
%Let $R = R_0, R_1, \ldots, R_k$ be a sequence of algebras with discrete rank $1$ valuations $v_i: R_i \setminus\{0\} \to Q$ for $0 \leq i < k$ such that $R_{i+1} = gr_{v_i}(R_i)$. Additionally, assume that $v_i$ is chosen to be homogeneous with respect to the $\Q^{i-1}$ grading on $R_i \cong gr_{v_{i-1}}(\cdots gr_{v_0}(R) \cdots)$, but is not a grading function.  Then there is a rank $k$ valuation $v: R\setminus\{0\} \to \Q^k$ with $gr_{v}(R) \cong R_k$.
%\end{theorem}
%\begin{proof}
%This follows from repeated application of Propositions \ref{concatenationgrading} and \ref{concatenationrank}. 
%\end{proof}

%Now let us consider the situation in Section \ref{subsec-degen-stages}. %Namely, we have algebras $R^(i)$. (but not positively graded).
%\begin{corollary}
%\end{corollary}

\begin{proof}[Proof of Theorem \ref{th-main}]
Following the proof of Lemma \ref{lem-degen-in-stages}, we can construct a sequence of algebras $R = R_0, R_1, \ldots, R_{d-1}$, and a sequence of valuations $v_i: R\setminus \{0\} \to \Z$ such that $\gr_{v_i}(R_i) = R_{i+1}$.   Here $v_0$ is the grading function given by the positive grading on $R$.  Each $v_i$ is constructed by the symbolic normal cone method, and is homogeneous with respect to a $\mathbb{G}_m^i$-action on $R_i$. Inductively, the $\mathbb{G}_m^i$-action on $R_i = \gr_{v_i}(R_{i-1})$ is constructed from the $\mathbb{G}_m^{i-1}$-action on $R_{i-1}$ and the $\mathbb{G}_m$-action on $\gr_{v_i}(R_{i-1})$ coming from its natural $\Z_{\geq 0}$-grading. The construction in the proof of Lemma \ref{lem-degen-in-stages} also implies that $v_i$ is not a grading function with respect to the $\mathbb{G}_m^i$-action. Repeated application of Proposition \ref{concatenationrank} now implies that the concatenation $\mathfrak{v} = v_0 \circ \ldots \circ v_{d-1}: R \setminus \{0\} \to \Z_{\geq 0} \times \Z^{d-1}$ has rank $d$. Furthermore, by Proposition \ref{prop-gr-concat} we have $\gr_{\mathfrak{v}}(R) \cong R_{d-1}$. Since $R_{d-1}$ is finitely generated by construction, it follows that $\mathfrak{v}$ has a finite Khovanskii basis.  
\end{proof}

We finish this section by showing that any valuation, whose value group is $\Q^r$ (for some $r >0$) and ordered lexicographically, can be obtained as a concatenation of rank $1$ valuations in the fashion described above.  This is a consequence of the following proposition. 

\begin{prop}\label{converse}
Let $\Q^r$ and $\Q^s$ be equipped with group orderings and order $\Q^r \times \Q^s$ lexicographically such that the $\Q^r$ component is preferred to the $\Q^s$ component. Let $\mathfrak{v}: R \setminus \{0\} \to \Q^r \times \Q^s$ be a valuation. Then there are valuations $v: R \setminus\{0\} \to \Q^r$ and $w: \gr_v(R) \setminus \{0\} \to \Q^s$ such that $w$ is homogeneous with respect to the grading on $\gr_v(R)$, and $\mathfrak{v} = v \circ w$. 
\end{prop}

\begin{proof}
First observe that the axioms of valuations implies that the projection $\pi_1 \circ \mathfrak{v} = v$ onto the first $r$ components is itself a valuation.  For $f \in R$ let $\mathfrak{v}(f) = (a, b)$, we claim that $b$ only depends on the equivalence class $\bar{f} \in gr_v(R)$.  If $v(f) = v(f') < v(f - Cf')$ for some $C \in k$, then in turn we must have $\mathfrak{v}(f - Cf') > \mathfrak{v}(f), \mathfrak{v}(f')$, but this implies that $\mathfrak{v}(f) = \mathfrak{v}(f')$.   Therefore, the projection $\pi_2\circ \mathfrak{v} = w$ is a well-defined function on homogeneous elements of $\gr_v(R)$; we extend this function to all of $\gr_v(R)$ using the $\min$ convention.  It is straightforward to verify that this function is super-additive.  We take two inhomogeneous elements $\sum \bar{f}$ and $\sum \bar{h}$.   We have $w(\bar{f}\bar{h}) = w(\bar{fh}) = w(\bar{f}) + w(\bar{h})$.   If we fix a grading degree $a$ and consider the sum $\sum \bar{f}_i\bar{h}_i$ of all products with degree $a$, we must have $w(\sum \bar{f}_i\bar{h}_i) \geq \min \{ w(\bar{f}_i) + w\bar{h}_i)\}$.  It follows that without loss of generality we may assume that all the elements in $\sum \bar{f}$ have the same $w$-value $b_1$, and likewise it can be assumed that all of the summands of $\sum \bar{h}$ have $w$-value $b_2$.   Now we consider the Newton polytopes of $\sum \bar{f}$ and $\sum \bar{h}$ in $\Q^r$.  Let $\bar{f}_e$ and $\bar{h}_e$ be components respectively at extremal vertices of these polytopes.  Then then $\bar{f}_e\bar{h}_e$ define an extremal vertex of the product $(\sum \bar{f})(\sum\bar{h})$, and the value $w(\bar{f}_e)\bar{h}_e)$ must be $b_1 + b_2 = w(\sum \bar{f}) + w(\sum \bar{h})$. 
\end{proof}

%\begin{remark}
%We translate the constructions in this section into the language of {\it weight valuations} used in \cite{Kaveh-Manon-NOK}.  Recall that if $R$ is positively graded and $\mathfrak{v}: R \setminus \{0\} \to \Z^d$ is homogeneous with finite Khovanskii basis $\mathcal{B} \subset R$, the $\mathfrak{v}$ is the {\it weight valuation} $\mathfrak{v}_M$, where $M$ is the matrix with columns $\{\mathfrak{v}(b_i)$ for $b_i \in \mathcal{B}$ (see \cite[Section 4]{Kaveh-Manon-NOK} and \cite[Lemma 3]{Kaveh-Manon-NOK}). 

%Following the proof of Theorem \ref{th-main}, Let $\bar{\mathcal{B}} \subset R_d$ be a finite generating set, and let $\mathcal{B} \subset R$ be a lift of this set, then $\mathcal{B}$ is a finite Khovanskii basis of $\mathfrak{v}$.   Now let $\mathcal{B}_i \subset R_i$ be a lift of $\bar{\mathcal{B}}$.  The valuation $v_i$ has Khovanskii basis $\mathcal{B}_i$, so it follows that $v_i = v_{u_i}$, where $u_i = (v_i(b_1), \ldots v_i(b_n))$.  By definition, the concatenation $\mathfrak{v} = v_0 \circ \cdots \circ v_d$ must be the weight valuation associated to $M = [\ldots u_i \ldots ]^t$. 
%\end{remark}

\section{Appendix: A lemma on finite generation of symbolic Rees algebra}     \label{sec-appendix}
In this appendix we prove the key lemma used in the proof of degeneration in stages (Lemma \ref{lem-degen-in-stages}).
\begin{lemma}  \label{lem-main}
Let $R$ be a finitely generated $\k$-domain. Suppose $0 \neq x \in R$ is such that:
\begin{itemize}
\item[(i)] The radical $\p = \sqrt{(x)}$ is a prime ideal.
\item[(ii)] The localization $R_\p$ is a discrete valuation ring.
\end{itemize}
Then the following hold:
\begin{itemize}
\item[(1)] The symbolic Rees algebra $\mathcal{A}_{\p^{(*)}}(R) = \bigoplus_k \p^{(k)}$ is a finitely generated $\k$-algebra.
\end{itemize}
Consider the homomorphism $\psi: (R/\p)[u] \to \gr_{\p^{(*)}} R$, where $(R/\p)[u]$ is the polynomial ring in one indeterminate $u$ over $R/\p$, given by $\psi(r) = r$, for $r \in R/\p$ and $\psi(u) = \bar{x}$, the image of $x$ in $\gr_{\p^{(*)}} R$.
\begin{itemize}
\item[(2)] The homomorphism $\psi$ is injective.
\item[(3)] $\gr_{\p^{(*)}} R$ regarded as a module over $(R/\p)[u]$ via $\psi$ is a finite module.
\item[(4)] If $R$ is moreover $\Z^r$-graded and $x$ is homogeneous with respect to this $\Z^r$-grading then the homomorphism $\psi$ is grading preserving. 
\item[(5)] If $R$ is positively graded and $x$ is homogeneous of degree $1$, then the embedding $\psi: (R/\p)[u] \hookrightarrow \gr_{\p^{(*)}} R$ induces a finite morphism $\phi: \Proj(\gr_{\p^(*)} R) \to  \Proj((R/\p)[u])$. 
\end{itemize}
\end{lemma}
We will need the following lemma.
\begin{lemma}   \label{lem-Veronese-subring}
Let $A = \bigoplus_{i \geq 0} A_i$ be a positively graded integral domain. If for some integer $m>0$, the $m$-th Veronese subring $A^{[m]} = \bigoplus_{i \geq 0} A_{mi}$ is a Noetherian ring, then $A$ is finite over $A^{[m]}$.
\end{lemma}
\begin{proof}
Let us write $A = \bigoplus_{j=0}^{m-1} M_j$ where $M_j = \bigoplus_{k = 0}^{\infty} A_{km + j}.$ If $M_j \ne \{0\}$, then it contains a nonzero homogeneous element $r$. Then the multiplication by $r^{m-1}$map $M_j \to A^{[m]}$ is well-defined and injective since $A$ is an integral domain. It is also a homomorphism of $A^{[m]}$-modules. It follows that $A$ is a finitely generated $A^{[m]}$-module as required.
\end{proof}

\begin{proof}[Proof of Lemma \ref{lem-main}] 
(1) Let $n>0$ be such that $x \in \p^{(n)} \setminus \p^{(n+1)}$. First we show that for all $k>0$ we have: 
$$\p^{(kn)} = \overline{(x^k)},$$ where $\overline{(x^k)}$ denotes the integral closure of the principal ideal $(x^k)$ in the field of fractions of $R$. Firstly, $\sqrt{\overline{(x^k)}} = \sqrt{(x^k)} = \p$. Thus $\p$ is an associated prime of $\overline{(x^k)}$. By  \cite[Corollary 5.4.2]{HS} we know that $\overline{(x^k)}$ is unmixed. It follows that in its primary decomposition only primary ideals corresponding to minimal primes appear. Thus $\overline{(x^k)} = \q$ where $\q$ is a $\p$-primary ideal. Moreover, $\q = \overline{(x^k)} A_{\p} \cap A$. Now $A_{\p}$ is a DVR so there exists $n > 0$ such that 
$$\overline{(x^k)} A_{\p} = (x^k)A_\p = \p^{kn} A_{\p}.$$
Here we are using the fact that any principal ideal in a normal ring is integrally closed. Thus $$\q = \p^{kn}A_{\p} \cap A = \p^{(kn)}$$ which proves the claim. We note that $S' = \bigoplus_{k} \overline{(x^k)} = \overline{\bigoplus_{k} (x^k)}$, the integral closure of the ring $S = \bigoplus_{k} (x^k)$ in $\bigoplus_k R$ (see \cite[Proposition 5.2.1]{HS}). By the finiteness of integral closure we know that $S'$ is finite over $S$. Now, the finite generation of $\bigoplus_{k} \p^{(k)}$ follows from Lemma \ref{lem-Veronese-subring}.

(2) By definition $\psi$ preserves the natural $\Z_{\geq 0}$-gradings on $(R/\p)[u]$ and $\gr_{\p^{(*)}} R$. So $\ker \psi$ is a homogeneous ideal in $(R/\p)[u]$. Let $(a+\p)u^k \in \ker \psi$, then $\bar{a}(\bar{x})^k = 0$ where $\bar{a}$, $\bar{x}$ are the images of $a$, $x$ in $\gr_{\p^{(*)}} R$ respectively. But since $\gr_{\p^{(*)}}R$ is a domain this is not possible unless $\ker \psi = \{0\}$. 

(3) Recall from proof of (1) that $S' = \bigoplus_{k} \overline{(x^k)} = \bigoplus_k \p^{(kn)}$ is finite over $S = \bigoplus_k (x^k)$. Let $\{f_1, \ldots, f_s\} \subset S'$ be a finite set of homogeneous $S$-module generators with $f_i \in \overline{(x^{k_i})}$. For $k>0$ take $f \in \overline{(x^k)} = \p^{(kn)}$ such that $f \notin \p^{(kn+1)}$, that is, $0 \neq \bar{f} \in \p^{(kn)} / \p^{(kn+1)}$. There exists $r_i \in R$ such that $$f = \sum_i r_ix^{k-k_i} f_i.$$
We recall that $x \in \p^{(n)}$. Going modulo $\p^{(kn+1)}$ we have: 
\begin{equation}  \label{equ-f}
f = \sum_j r_j x^{k-k_j} f_j \quad (\textup{mod}~ \p^{(kn+1)}),
\end{equation}
where the sum is over all $j$ such that $r_j \notin \p$ and $f_j \notin \p^{(k_jn+1)}$. Now by definition $\bar{f}$ (respectively $\bar{f}_j$) is the image of $f$ in $\p^{(kn)} / \p^{(kn+1)}$ (respectively $\p^{(k_jn)} / \p^{(k_jn+1)}$). It follows from \eqref{equ-f} that $\bar{f}$ is a linear combination of the $\bar{f}_j$ with coefficients $\psi(r_ju^{k-k_j}) = (r_j + \p)\bar{x}^{k-k_j}$ respectively. This shows that the $\bar{f}_i$ are a finite set of module generators for $\bigoplus_k \p^{(kn)} / \p^{(kn+1)}$ over $(R/\p)[u]$. As before, the claim follows from Lemma \ref{lem-Veronese-subring}.  

(4) It follows from the construction of $\psi$. (5) We show the following more general fact. Let $A$, $B$ be positively graded domains and $\psi: A \to B$ an injective finite homomorphism. Then $\psi$ induces a morphism $\phi: \Proj(B) \to \Proj(A)$. For this, consider the morphism $\psi^*: \Spec(B) \to \Spec(A)$. Finiteness of $\psi$ implies that the inverse image (under $\psi^*$) of any point is a finite set and hence $\psi^*$ cannot send a $\G_m$-invariant subvariety of $\Spec(B)$ of positive dimension to the vertex $V(A_+)$ where $A_+$ is the irrelevant ideal,  generated by all the elements in $A$ of positive degree. It follows that the inverse image of the vertex of $\Spec(A)$ is the vertex of $\Spec(B)$ and hence $\psi^*$ induces a morphism $\phi: \Proj(B) \to \Proj(A)$. The finiteness of $\phi$ follows from the finiteness of $\psi$.
\end{proof}

\end{document}